\documentclass[oneside, 12pt]{amsart}

\usepackage{amscd, amssymb, amsmath, mathrsfs}
\usepackage[english]{babel}
\usepackage{booktabs}
\usepackage{tikz-cd}
\usepackage[colorlinks, linktocpage, citecolor = blue, linkcolor = blue]{hyperref}

\newcommand\mylabel[1]{\label{#1}}

\newtheorem{theorem}{Theorem}[section]
\newtheorem*{maintheorem}{Theorem}
\newtheorem{lemma}[theorem]{Lemma}
\newtheorem{proposition}[theorem]{Proposition}
\newtheorem{corollary}[theorem]{Corollary}

\theoremstyle{definition}
\newtheorem{definition}[theorem]{Definition}

\newtheorem*{acknowledgement}{Acknowledgement}

\theoremstyle{remark}

\newcommand{\NN}	{\mathbb{N}}
\newcommand{\ZZ}	{\mathbb{Z}}
\newcommand{\QQ}	{\mathbb{Q}}

\newcommand{\FF}	{\mathbb{F}}

\newcommand{\PP}	{\mathbb{P}}

\newcommand{\GG}	{\mathbb{G}}

\newcommand{\ideala}    {\mathfrak{a}}

\newcommand  {\shE}     {\mathscr{E}}
\newcommand  {\shF}     {\mathscr{F}}

\newcommand  {\shI}     {\mathscr{I}}

\newcommand  {\shM}     {\mathscr{M}}

\newcommand  {\shN}     {\mathscr{N}}
\newcommand  {\shL}     {\mathscr{L}}

\newcommand  {\shS}     {\mathscr{S}}
\newcommand  {\shT}     {\mathscr{T}}

\newcommand  {\calA}     {\mathcal{A}}

\newcommand  {\calC}     {\mathcal{C}}

\newcommand  {\foA}     {\mathfrak{A}}

\newcommand  {\foG}     {\mathfrak{G}}
\newcommand  {\foH}     {\mathfrak{H}}

\newcommand  {\foL}     {\mathfrak{L}}

\newcommand  {\foX}     {\mathfrak{X}}
\newcommand  {\foY}     {\mathfrak{Y}}

\newcommand{\Aff}	{\text{\rm Aff}}

\newcommand{\alg}	{\text{\rm alg}}

\newcommand{\Aut}	{\operatorname{Aut}}

\newcommand{\can}	{\text{\rm can}}

\newcommand{\cl}	{\operatorname{cl}}

\newcommand{\Cokernel}	{\operatorname{Coker}}

\newcommand{\Ext}	{\operatorname{Ext}}

\newcommand{\Frac}	{\operatorname{Frac}}

\newcommand{\GL}	{\operatorname{GL}}

\newcommand{\Hom}	{\operatorname{Hom}}

\newcommand{\id}	{{\operatorname{id}}}
\newcommand{\Image}	{\operatorname{Im}}

\newcommand{\II}  	{\text{\rm II}}

\newcommand{\Kernel} 	{\operatorname{Ker}}

\newcommand{\invlim}	{\varprojlim}

\newcommand{\lra}	{\longrightarrow}

\newcommand  {\maxid}   {\mathfrak{m}}

\newcommand  {\mult}    {\text{\rm mult}}

\newcommand  {\Num}     {\operatorname{Num}}

\newcommand  {\Nil}     {\operatorname{Nil}}

\newcommand  {\primid}  {\mathfrak{p}}
\renewcommand{\O}       {\mathscr{O}}

\newcommand  {\ord}     {\operatorname{ord}}

\newcommand  {\Pic}     {\operatorname{Pic}}

\newcommand  {\pr}      {\operatorname{pr}}

\newcommand  {\quadand} {\quad\text{and}\quad}

\newcommand  {\ra}      {\rightarrow}

\newcommand  {\rank}    {\operatorname{rank}}
\newcommand  {\red}     {{\operatorname{red}}}

\newcommand  {\sep}     {{\operatorname{sep}}}

\newcommand  {\Sing}    {\operatorname{Sing}}

\newcommand  {\Spec}    {\operatorname{Spec}}

\newcommand  {\Spf}     {\operatorname{Spf}}

\newcommand{\uHom}{\underline{\operatorname{Hom}}}

\newcommand{\mydate}{
\number\day\space
\ifcase\month \or January\or February\or March\or April\or May\or June\or July\or August\or September\or October\or November\or December\fi 
\space\number\year
}

\newcommand{\Cp}{\operatorname{Cp}}
\newcommand{\shCp}{\mathscr{C}\!\text{\textit{p}}}
\newcommand{\shSc}{\mathscr{S}\!\text{\textit{c}}}
\newcommand{\shRel}{\mathscr{R}\!\text{\textit{el}}}

\begin{document}

\title[The Deligne--Illusie Theorem]
      {The Deligne--Illusie Theorem and exceptional Enriques surfaces}

\author[Stefan Schr\"oer]{Stefan Schr\"oer}
\address{Mathematisches Institut, Heinrich-Heine-Universit\"at,
40204 D\"usseldorf, Germany}
\curraddr{}
\email{schroeer@math.uni-duesseldorf.de}

%formal methods and deformations, K3 surfaces and Enriques surfaces, vector bundles, group schemes, elliptic surfaces, exact and abelian categories

%\dedicatory{Preliminary version, \mydate}
\dedicatory{Revised version, 4 January 2021}

\begin{abstract}
Building on the results of Deligne and Illusie on liftings to
truncated Witt vectors, we give a criterion for non-liftability 
that involves only the dimension of certain cohomology groups of  
vector bundles arising from the Frobenius pushforward of the de Rham complex.
Using vector bundle methods, we apply this to show that  exceptional Enriques surfaces, a class introduced by
Ekedahl and Shepherd-Barron, do not lift to truncated Witt vectors, yet the base of the miniversal
formal deformation over the Witt vectors is regular.
Using the classification of Bombieri and Mumford, we also show  that bielliptic
surfaces arising from a quotient by a  unipotent group scheme of order $p$ do not
lift to the ring of Witt vectors.
These results hinge on some observations in homological algebra that relates splittings in
derived categories to Yoneda extensions and certain diagram completions.
\end{abstract}

\maketitle
\tableofcontents

%===========================================================
\section*{Introduction}
\mylabel{Introduction}

Let $k$ be a perfect field of characteristic $p>0$, and $Y$ be a smooth proper $k$-scheme.
Often it is a challenging question  whether or not the scheme $Y$ lifts to the ring of
Witt vectors $W$, or even its truncation  
$W_2=W/p^2W$.
According to a famous result  of Deligne and Illusie \cite{Deligne; Illusie 1987}, the existence of such   $W_2$-liftings  
implies that   the Hodge--de Rham spectral sequence 
$E_1^{rs}=H^s(Y,\Omega_Y^r)\Longrightarrow H^{r+s}(Y,\Omega_Y^\bullet)$ degenerates at the $E_1$-page, provided $\dim(Y)\leq p$.
This result is  used to show that   schemes with ``exotic''   Hodge cohomology frequently do  not admit such   lifts.
Consequently, the base $\Spf(A)$ of the miniversal formal deformation is not formally smooth   over the ring $W$.
Note that it could still be given by a regular local ring, for example $A=W[[U,V]]/(UV-p)$.
Also note that failure of lifting to $W_2$ occurs  in surprisingly simple situations, even for smooth models of inseparable covers of the projective plane 
(\cite{Liedtke; Strano 2014}, Theorem 3.4).

Actually, Deligne and Illusie identified the \emph{gerbe of liftings} of the scheme $Y'$
to the ring $W_2$   with the \emph{gerbe of splittings} for the one-term complex
$F_*\O_Y\stackrel{d}{\ra} Z\Omega^1_Y$.
Here $Y'=Y\otimes_kk$ is the base-change with respect to the Frobenius map $\lambda\mapsto\lambda^p$,
and  $Z\Omega^1_Y$ is the sheaf of 1-cocycles
in the   push-forward $F_*(\Omega_Y^\bullet)$ of the de Rham complex with respect to the relative Frobenius $F:Y\ra Y'$.
To my best knowledge, this amazing result was \emph{never  used directly} to show that certain schemes  do not lift to the ring $W_2$.
The main goal of this paper is to show that such arguments are indeed feasible.
For this, we establish \emph{general numerical criteria}  that ensure that the gerbes in question have no global objects,
and apply this to certain  \emph{  Enriques surfaces}  and     \emph{bielliptic surfaces}.

This hinges on some general results in homological algebra, which ensure among other things that the above gerbes
admit a global object if and only if the Yoneda class of the four-term exact sequence
$$
0\lra \O_{Y'}\lra F_*(\Omega_Y)\stackrel{d}{\lra} Z\Omega^1_Y\lra \Omega^1_{Y'}\lra 0
$$
in $\Ext^2(\Omega^1_{Y'},\O_{Y'})$ vanishes. This sequence is obtained by splicing two short exact sequences 
$$
0\ra  \O_{Y'}\ra F_*(\Omega^1_Y)\ra B\Omega^1_Y\ra 0\quadand
0\ra  B\Omega^1_Y\ra Z(\Omega^1_Y)\ra  \Omega^1_{Y'}\ra 0.
$$
If the former splits, one says that $Y$ is \emph{Frobenius-split}. This notation was introduced by Mehta and Ramanathan \cite{Mehta; Ramanathan 1985},
and has numerous striking  consequences for the cohomology of  sheaves.
We refer to the monograph of Brion and Kumar \cite{Brion; Kumar 2005} for a highly readable account.
Let us say that $Y$ is \emph{Cartier-split} if the
second short exact sequence splits. As Srinivas \cite{{Srinivas 1990}} and Yobuko \cite{Yobuko 2019} observed, this means  
that the scheme $Y$ and also  the morphism $F:Y\ra Y'$ admits a lifting to the ring $W_2$.

Here we are interested in a much weaker and more flexible version:
We say that the scheme $Y$ is \emph{pre-Cartier-split} if the 
connecting map $\Hom(B\Omega^1_Y,\O_{Y'})\ra \Ext^1(\Omega^1_{Y'},\O_{Y'})$ coming from the second short exact sequence is the zero map.
Our first main result is the following general numerical criterion:

\begin{maintheorem}
(See   \ref{numerical criterion}.)
Suppose $Y$ is pre-Cartier-split but not Frobenius-split, and satisfies
$h^1(\Theta_Y)\geq h^1(\uHom(Z\Omega^1_Y,\O_{Y'}))$.
Then   the scheme $Y$ does not lift to the ring of truncated Witt vectors $W_2$.
\end{maintheorem}

Under the condition $c_1=0$ and $\dim(Y)=2$ this    simplifies further.
One gets the following   version, in which the differentials of the de Rham complex are eliminated:

\begin{maintheorem}
(See \ref{non-liftable surface}.)
Suppose $Y$ is a surface, that the dualizing sheaf $\omega_Y$ has order  $p\geq 2$ in the Picard group,
and that $h^1(\Theta_Y) \geq h^1(\uHom(F_*\Omega^1_Y,\O_{Y'}))$. 
Then the scheme $Y$ does not lift to the ring $W_2$.
\end{maintheorem}

This is amenable to vector bundle techniques, and the main idea is to put the cotangent sheaf into
a short exact sequence $0\ra\shL\otimes\omega_Y\ra \Omega^1_Y\ra \shI\shL^\vee\ra 0$,
and write   $\uHom(F_*\Omega^1_Y,\O_{Y'})$ on $Y'$
as  the Frobenius image of some   sheaf on $Y$, in order to compute cohomology.
In characteristic $p=2$, we    formulate in Theorem \ref{special surface no lift} four elementary conditions
concerning the invertible sheaves $\omega_Y$  and   $\shL$ that  ensure that the surface $Y$ does not lift to the ring $W_2$.
 
We then apply our results to \emph{exceptional Enriques surface} $Y$.
These surfaces were introduced and studied by Ekedahl and Shepherd-Barron \cite{Ekedahl; Shepherd-Barron 2004}.
Despite their highly unusual geometry, which was already considered by Cossec and Dolgachev \cite{Cossec; Dolgachev 1989},
they admit rather concrete descriptions in terms of   equations found by Salomonsson \cite{Salomonsson 2003}.
Actually, one should treat them together with the \emph{supersingular Enriques surfaces}, because both
have $h^0(\Theta_Y)=h^2(\Theta_Y)=1$ and $h^1(\Theta_Y)=12$, which makes their deformation theory seemingly complicated.
Ekedahl, Hyland and Shepherd-Barron  \cite{Ekedahl; Hyland; Shepherd-Barron 2012} showed that
the base of the miniversal formal deformation   of a supersingular Enriques surface
is the formal spectrum of $A=W[[T_1,\ldots,T_{12}]]/(2-FG)$ where $F,G$ lie in the maximal ideal. Before Theorem 4.5, they ask
to clarify the miniversal deformation in the exceptional case.  Extending their results, we  obtain: 

\begin{maintheorem}
(See \ref{versal deformation enriques}.)
The base of the miniversal deformation of an Enriques surface $Y$ in characteristic two
is given by a complete local noetherian ring  $A$  that is regular, $W$-flat and 
of Krull dimension eleven or twelve. Moreover, the following are equivalent:
\begin{enumerate}
\item The Enriques surface $Y$ is exceptional or supersingular.
\item The scheme $Y$ does not lift  to the ring $W_2$.
\item The absolute ramification index is $e(A)\geq 2$.
\item The dimension is $\dim(A)=12$.
\end{enumerate}
\end{maintheorem}

This is in striking contrast to general results of Liedtke \cite{Liedtke 2015}, who showed that
\emph{normal} Enriques surfaces having a Cossec--Verra polarization are unobstructed.
The non-liftability of our smooth $Y$  thus  must
be caused by the necessity of 
base-change needed in Artin's simultaneous resolution   \cite{Artin 1974} of the 
singularities in some normal models of $Y$, all of which are rational double points. This   was further
elucidated by Shepherd-Barron \cite{Shepherd-Barron 2017}.
% The non-liftability of $Y$  thus  must
% be caused by properties of rational double points in some normal models, coming from the
% base-change needed for Artin's simultaneous resolution of singularities \cite{Artin 1974}. This   was further
% elucidated by Shepherd-Barron \cite{Shepherd-Barron 2017}.
One should compare the above result with the situation in characteristic zero:
Then  the $T^1$-lifting Theorem ensures that the base of the miniversal deformation for
smooth schemes with $c_1=0$  is given by a regular ring  (confer \cite{Bogomolov 1978},  \cite{Todorov 1980}, \cite{Tian 1987},
\cite{Kawamata 1992},  \cite{Ran 1992},  \cite{Schroeer 2003}, \cite{Ekedahl; Shepherd-Barron 2003}).

Finally, we apply our results to \emph{bielliptic surfaces} $Y$, which were
classified by Bombieri and Mumford \cite{Bombieri; Mumford 1977}, \cite{Bombieri; Mumford 1976}.
Their deformation theory was studied by Partsch \cite{Partsch 2013},
when both genus-one fibrations are elliptic.
Here we examine the case that the surface is of the form $Y=(E\times C)/G$, where
$C$ is the rational cuspidal curve,    
 $G$ is a finite group scheme, and the characteristic is $p=2$.
We shall see that if $G=\alpha_2$, these surfaces do not lift to the ring
of Witt vectors $W$. We also describe the tangent and cotangent sheaves
and their cohomology. The assertion depends on   general results about proper
group schemes and Picard schemes, in particular:

\begin{maintheorem}
(see   \ref{picard no lift}.)
Set $G=\Pic^0_{Y/k}$. Suppose the local  group scheme $L=G/G_\red$ contains
some $\alpha_{p^n}$, $n\geq 1$ as a direct summand,
and that  the first Betti number satisfies  $b_1=2(h^1(\O_Y)-h^2(\O_Y))$. Then the scheme $Y$ does not lift to the ring $W$.
\end{maintheorem}

All the above result hinge on certain general results from homological algebra,
in which we give a new interpretation of the gerbe of splittings for a two-term complex $f:M\ra N$
in some abelian category $\calA$. Here we introduce the notion of \emph{diagram completions},
which consists of an object $E$, together with two morphisms $h:M\ra E$ and $g:E\ra N$
making 
\begin{equation} 
\begin{CD}
M	@>\pr>>	B\\
@VhVV		@VViV\\
E	@>>g>	N
\end{CD}
\end{equation}
both cartesian and cocartesian. If the objects of the abelian
category are abelian sheaves on some space or site, and the cohomology sheaves $H^i$
for the one-term complex $M\stackrel{f}{\ra} N$ are locally free of
finite rank, we get the following
result,  which seems to be of independent interest:

\begin{maintheorem}
(See  \ref{same gerbe}.)
The gerbe   of diagram completions and the gerbe  of splittings for $f:M\ra N$ have the same class in 
the   group
$$
H^2(X,\uHom_{\O_X}(H^1,H^0))=\Ext^2(H^1,H^0).
$$
Moreover, either  of them admits a global object if and only if the Yoneda class  
of the exact sequence $0\ra H^0\ra M\stackrel{f}{\ra} N\ra H^1\ra 0$   vanishes.
\end{maintheorem}

\medskip
This paper is organized as follows:
In Section \ref{Numerical criteria} we introduce the notion of pre-Cartier split schemes,
and give  our general numerical criterion against liftings to the ring $W_2$.
In Section \ref{Condition c1} this is examined under the additional conditions $c_1=0$
and $\dim(Y)=2$.
Section \ref{Base miniversal} contains a discussion for the regularity properties
of the base of the miniversal deformation in mixed characteristics.
In Section \ref{Group schemes} we consider obstructions to liftings arising from
the theory of group schemes and Picard groups.
Section \ref{Vector bundle} contains some computations with vector bundles and Chern classes
on surfaces.
These are used in Section \ref{Algebraic surfaces}, to prove our main result on 
non-liftability of surfaces in characteristic two.
Section \ref{Enriques surfaces} and \ref{Bielliptic surfaces} contain the applications to Enriques and bielliptic surfaces.
The final Section \ref{Homological algebra} deals with necessary homological algebra in a general abstract setting.

\begin{acknowledgement}
I wish  to thank Luc Illusie for   valuable discussions, and the referees for thorough reading and helpful comments.
This research was conducted in the framework of the   research training group
\emph{GRK 2240: Algebro-geometric Methods in Algebra, Arithmetic and Topology}, which is funded
by the DFG. 
\end{acknowledgement}

%===========================================================
\section{Numerical criteria against first-order liftings}
\mylabel{Numerical criteria}

Let $k$ be a perfect field of characteristic $p>0$, and $Y$ be a smooth proper $k$-scheme 
of dimension $n=\dim(Y)$. Furthermore assume that $h^0(\O_Y)=1$.
To simplify notation, we write $\Omega^1_Y=\Omega^1_{Y/k}$ for the cotangent sheaf,
$\Theta_Y=\uHom(\Omega^1_Y,\O_Y)$ for the tangent sheaf,   $\Omega^i_Y=\Lambda^i(\Omega^1_Y)$ for the sheaves of differential forms,
and $Y'=Y\otimes_kk$ for the base-change with respect to the Frobenius map $\lambda\ra\lambda^p$.
The commutative diagram  
\begin{equation}
\label{relative frobenius}
\begin{tikzcd}[column sep = small]
Y\ar[rr, "F"]\ar[rrr,bend left=20, "F_Y"]\ar[dr]	& 			& Y'\ar[r]\ar[dl]	& Y\ar[d]\\
					& \Spec(k)\ar[rr, "F_k"']		&& \Spec(k),
\end{tikzcd}
\end{equation}
where $F_Y$ and $F_k$ are absolute Frobenius morphisms, 
defines the \emph{relative Frobenius morphism} $F :Y\ra Y'$, which is a finite flat universal homeomorphism of 
degree $\deg(Y/Y')=p^n$.
Note  that for each skyscraper sheaf
$\shT$ on $Y$  the Frobenius pushforward $F_*(\shT)$ is a skyscraper sheaf on $Y'$ with $ h^0(F_*\shT)=h^0(\shT)$,
and that for each locally free $\O_Y$-module $\shE$ of rank $r\geq 0$ the Frobenius
pushforward $\shE'=F_*(\shE)$ is a locally free $\O_{Y'}$-module of rank $r'=rp^n$

The latter applies in particular to the $\shE=\Omega_Y^i$. Moreover, the $k$-linear differentials
in the de Rham complex $\Omega^\bullet_Y$ become $\O_{Y'}$-linear maps in the resulting cochain complex $F_*(\Omega^\bullet_Y)$.
We write $B\Omega_Y^i$ and $Z\Omega_Y^i$ for   coboundaries and cocycles, viewed as coherent sheaves on $Y'$.
These are actually locally free, and the \emph{inverse Cartier operator} $d(f\otimes1)\mapsto f^{p-1}df$ gives an identification
$\Omega^i_{Y'}=Z\Omega^i_Y/B\Omega^i_Y$. In particular, we have a commutative diagram
\begin{equation}
\label{four-term sequence}
\begin{tikzcd}[row sep=tiny, column sep=small ]
		&			&  		 0\arrow[dr]			& 		& 0 \\		
		&			&						& B\Omega^1_Y\arrow[dr]\arrow[ur]\\
0\arrow[r]	& \O_{Y'}\arrow[r]	& F_*\O_Y\arrow[ru]\arrow[rr, "d"']		& 		& Z\Omega^1_Y\arrow[r]	&\Omega^1_{Y'}\arrow[r]	& 0,\\
\end{tikzcd}
\end{equation}
where the four-term horizontal sequence is exact, and arises  via splicing of the two short exact sequences with kinks.

We now consider the second cohomology group
$$
H^2(Y',\Theta_{Y'})=\Ext^2(\Omega_{Y'}^1,\O_{Y'}), 
$$
whose elements can be regarded as equivalence classes of \emph{gerbes banded by $\Theta_{Y'}$}, or
in short $\Theta_{Y'}$-gerbes 
(\cite{Giraud 1971}, Chapter IV, Section 3.4). It contains the class of the 
\emph{gerbe  of splittings} $\shSc(F_*\O_Y\stackrel{d}{\ra} Z\Omega^1_Y)$ for  the two-term complex
$F_*\O_Y\stackrel{d}{\ra} Z\Omega^1_Y$, as defined in \cite{Deligne; Illusie 1987}, Section 3.
It also contains the class of the \emph{gerbe of liftings} $\shRel(Y',W_2)$ for the $k$-scheme $Y'$ to the ring $W_2$.
The objects of this gerbe are the 
proper flat morphisms $\foY'\ra\Spec(W_2)$, together with an identification $\foY'\otimes_{W_2}k=Y'$.
Here   $W_2=W/p^2W$ is the truncation of length two for the ring of Witt vectors $W$.
We refer to \cite{AC 8-9}, Chapter IV, \S1 for a comprehensive treatment of  Witt vectors.
Since the  Frobenius map $k\ra k$ is bijective, it induces an automorphism of the local ring $W$,
and the liftings of $Y'$ to $W_2$ correspond to  the liftings of $Y$.
According to \cite{Deligne; Illusie 1987}, Proposition 3.3 there is an equality
$$
\cl\shSc(F_*\O_Y\stackrel{d}{\ra} Z\Omega^1_Y) = \cl\shRel(Y',W_2)  
$$
of classes in the second cohomology group $H^2(Y',\Theta_{Y'})$. From the abstract situation treated
in Theorem \ref{same gerbe} below we get
another, completely different interpretation:

\begin{proposition}
\mylabel{insert cartesian}
There is a lifting of the scheme $Y$ to the ring $W_2$ if and only if the Yoneda class of the
four-term exact sequence in \eqref{four-term sequence}  vanishes  in the extension group $\Ext^2(\Omega^1_{Y'},\O_{Y'})$.
Up to isomorphism, to give such a lifting  amounts to factor the differential
$d:F_*\O_Y\ra Z\Omega^1_Y$ over some locally free $\O_Y'$-module  $\shE$ of rank $p^n+n$ such that the diagram
$$
\begin{CD}
F_*\O_Y	@>>>	B\Omega^1_Y\\
@VVV		@VVV\\
\shE	@>>>	Z\Omega^1_Y
\end{CD}
$$
becomes both cartesian and cocartesian.
\end{proposition}

In somewhat different form, this was already observed by De Clercq, Florence and Lucchini Arteche 
(\cite{De Clercq; Florence; Lucchini Arteche 2018}, Corollary 5.18).

The scheme $Y$ is called \emph{Frobenius-split} if the inclusion $\O_{Y'}\subset F_*(\O_Y)$ admits a retraction
\cite{Mehta; Ramanathan 1985}. In other words, the   extension class
$$
\cl(F_*\O_Y)\in\Ext^1(B\Omega^1_Y,\O_{Y'}) = H^1(Y',\uHom(B\Omega^1_Y,\O_{Y'}))
$$
vanishes. Let us say that $Y$ is \emph{Cartier-split} if the surjection $Z\Omega^1_Y\ra\Omega^1_{Y'}$ admits
a section, that is, the extension class
$\cl(Z\Omega_Y^1) \in \Ext^1(\Omega^1_{Y'},B\Omega^1_Y)$
is zero. Using that the \emph{Yoneda product} $\cl(F_*\O_Y) \ast \cl(Z\Omega^1_Y)$ is the \emph{Yoneda class} of the
four-term exact sequence in \eqref{four-term sequence}, we get with  Proposition \ref{insert cartesian}:

\begin{proposition}
\mylabel{frobenius and cartier}
If $Y$ is Frobenius-split or Cartier-split, then the scheme $Y$ lifts to the ring $W_2$.
\end{proposition}

We now introduce another condition that is much weaker than Cartier-split:
The short exact sequence to the right  in \eqref{four-term sequence} yields an exact sequence
\begin{equation}
\label{connecting map}
\Hom(Z\Omega^1_Y,\O_{Y'})\lra \Hom(B\Omega^1_Y,\O_{Y'})\stackrel{\partial}{\lra} \Ext^1(\Omega^1_{Y'},\O_{Y'}).
\end{equation}
The map on the right is the \emph{connecting map}, which may or may not vanish.

\begin{definition}
\mylabel{pre-cartier-split}
We say that $Y$ is \emph{pre-Cartier-split} if the connecting map $\partial$ 
in the exact sequence \eqref{connecting map} is the zero map.
\end{definition}

Saying that the scheme  $Y$ is Cartier-split  means that the 
short exact sequence $0\ra B\Omega_Y^1\ra Z\Omega^1_Y\ra\Omega^1_{Y'}\ra 0$ splits;
then the connecting map is a priori zero and $Y$ is also pre-Cartier-split.
The latter also holds if the group $\Hom(B\Omega^1_Y,\O_{Y'})$ vanishes.
A key observation for this paper is the following    criterion:

\begin{theorem}
\mylabel{numerical criterion}
Suppose $Y$ is pre-Cartier-split but not Frobenius-split, and satisfies
\begin{equation}
\label{inequality}
h^1(\Theta_Y)\geq h^1(\uHom(Z\Omega^1_Y,\O_{Y'})).
\end{equation}
Then this   inequality is an equality, and  the scheme $Y$ does not lift to the ring of truncated Witt vectors $W_2$.
\end{theorem}

\proof
Since $Y$ is not Frobenius-split, the extension class $\cl(F_*\O_Y)\in \Ext^1(B\Omega^1_Y,\O_{Y'})$
does not vanish. The   short exact sequence to the right in \eqref{four-term sequence} yields a long exact sequence
$$
\Ext^1(\Omega_{Y'}^1,\O_{Y'})\lra \Ext^1(Z\Omega^1_Y,\O_{Y'})\lra \Ext^1(B\Omega^1_Y,\O_{Y'})\lra\Ext^2(\Omega^1_{Y'},\O_{Y'}).
$$
The connecting map in \eqref{connecting map} vanishes by assumption, so the map on the left is injective.
It is actually bijective,  because by assumption the vector space dimension $h^1(\Theta_Y)$ of its domain
is at least as large as the dimension   of its range. In turn, the inequality \eqref{inequality} is an equality.
Moreover, the map in the middle of the above exact sequence is zero, so the extension class $\cl(F_*\O_Y)\neq 0$ is not in the image.
Hence its image in $\Ext^2(\Omega^1_{Y'},\O_{Y'})$ is non-zero.
This image equals the Yoneda class of the four-term sequence in  \eqref{four-term sequence}, by definition of the
connecting map for Yoneda extension groups.
According to Proposition \ref{insert cartesian}, the scheme $Y$ does not lift to 
the ring $W_2$. 
\qed

\medskip
This result reveals that under suitable assumptions, a mere bound on certain cohomology groups implies
non-existence of liftings.
Note that the scheme $Y$ is Frobenius split or Cartier split if and only if the respective property holds for the base-change 
to the algebraic closure $k^\alg$.

We are particularly interested in the situation where  $\omega_Y^{\otimes p}\simeq\O_Y$ but $\omega_Y\not\simeq\O_Y$.
Then the scheme $Y$   comes with a \emph{canonical covering}
$\epsilon:X\ra Y$, which is a torsor under the local group scheme $\mu_p=\uHom(\ZZ/p\ZZ,\GG_m)$.
This     explains why we choose
the symbol $Y$ for our smooth proper scheme. Note that the canonical covering is a 
finite flat universal homeomorphism of degree $\deg(X/Y)=p$, and that the scheme 
$X$ usually contains singularities, and easily  may become non-normal. 
For more details, see   for example   \cite{Schroeer 2017}, Section 4.

%===========================================================
\section{The condition \texorpdfstring{$c_1=0$}{c1=0}}
\mylabel{Condition c1}

We keep the assumptions of the previous section, such that $Y$ is a smooth proper $n$-dimensional  scheme
over some perfect field $k$ of characteristic $p>0$, with Frobenius pullback $Y'=Y\otimes_kk$.
The dualizing sheaf $\omega_Y=\Omega^n_Y$ and the relative dualizing
sheaf $\omega_{Y/Y'}$ will be of paramount importance. The latter is defined as a coherent sheaf on $Y$ by
the formula $F_*(\omega_{Y/Y'}) = \uHom(F_*\O_Y,\O_{Y'})$. Let us start with the following facts:

\begin{lemma}
\mylabel{dualizing sheaves}
The fibers of the relative Frobenius $F:Y\ra Y'$ are Gorenstein, and the
relative dualizing sheaf  
is given by  $\omega_{Y/Y'}=\omega_Y^{\otimes 1-p}$.
If $\omega_Y$ is $p$-torsion in the Picard group, we have 
$F_*(\omega_Y)=\underline{\Hom}(F_*\O_Y,\O_{Y'})$.
\end{lemma}

\proof 
Disregarding the structure morphisms to $\Spec(k)$, we first note that the projection $Y'=Y\otimes_kk\ra Y$ is an isomorphism of  schemes.
It thus suffices to verify the first statement for the absolute Frobenius $F_Y: Y\ra Y$.
Fix a   point $a\in Y$, and choose a regular system of parameters $f_1,\ldots,f_r\in \O_{Y,a}$.
The schematic fiber $Y_a=F_Y^{-1}(a)$ is the spectrum of the ring $R=\O_{Y,a}/(f_1^p,\ldots,f_r^p)$.
Passing to formal completions, we see that $R=\kappa(a)[[T_1,\ldots,T_r]]/(T_1^p,\ldots,T_r^p)$.
The socle of this local Artin ring is generated by $ \prod_{i=1}^rT_i^{p-1}$, whence the local ring
$R$ is Gorenstein.    

In turn, the relative dualizing sheaf is invertible, and satisfies the formula $\omega_Y=\omega_{Y/Y'}\otimes F^*(\omega_{Y'})$.
We have $\omega_{Y'}=\omega_Y\otimes_kk$, for the base-change with the Frobenius map $\lambda\mapsto \lambda^p$.
Since the absolute Frobenius induces multiplication-by-$p$ on the Picard group,
the diagram \eqref{relative frobenius} yields $F^*(\omega_{Y'})=\omega_Y^{\otimes p}$, and the
assertion on the relative dualizing sheaf follows.
Finally, if $\omega_Y$ is $p$-torsion, we get $\omega_{Y/Y'}=\omega_Y$, and the last statement
comes from the definition of   relative dualizing sheaves.
\qed
 
\medskip
Write $c_1= c_1(Y)=c_1(\Omega^1_Y)=c_1(\omega_Y)$ for the \emph{first Chern class},   say as an element in
the group $\Num(Y)$ of invertible sheaves modulo numerical equivalence. In other words, the condition $c_1=0$
means that  $(\omega_Y\cdot C)=0$ for all curves $C\subset Y$. 

\begin{lemma}
\mylabel{not one-dimensional}
Suppose $c_1=0$ holds. Then the vector space $\Hom(B\Omega^1_Y,\O_{Y'})$ is at most one-dimensional.
It is one-dimensional if and only if the scheme $Y$ is not Frobenius-split and the dualizing sheaf $\omega_Y$ is $(p-1)$-torsion 
in the Picard group. 
\end{lemma}

\proof
The   short exact sequence to the left in the diagram \eqref{four-term sequence} splits if and only if
the dual exact sequence  
$$
0\lra \uHom(B\Omega^1_Y,\O_{Y'})\lra F_*(\shL)\lra \O_{Y'}\lra 0
$$
splits, for the numerically trivial sheaf $\shL=\omega_{Y/Y'}=\omega_Y^{\otimes 1-p}$.
In turn, we get an inclusion $\Hom(B\Omega^1_Y,\O_{Y'})\subset H^0(Y,\shL)$.  Since $Y$ is integral,
we have $h^0(\shL)\leq 1$; equality holds if  and only if $\shL\simeq\O_Y$.
This shows that  $ \Hom(B\Omega^1_Y,\O_{Y'})$ is at most one-dimensional.

Suppose that it is one-dimensional. Then $h^0(\shL)\neq 0$, hence $h^0(\shL)=1$, and $\omega_Y$ is $(p-1)$-torsion.
Moreover, the extension does not split, because otherwise  the contribution of $h^0(\O_{Y'})=1$ yields
the contradiction $h^0(\shL)=2$.
Conversely, suppose that $\omega_{Y}^{\otimes 1-p}=\O_Y$ and that the extension does not split. 
In the resulting exact sequence
$$
0\lra \Hom(B\Omega^1_Y,\O_{Y'})\lra H^0(Y,\O_Y)\lra H^0(Y',\O_{Y'}),
$$
the unit section $1\in H^0(Y',\O_{Y'})$ is not in the image of the  map on the right.
So this map vanishes, and it follows that the term on the left is one-dimensional.
\qed

\begin{proposition}
\mylabel{pre-cartier not frobenius}
Suppose that the order of  $\omega_Y$  in the Picard group is $p^\nu$ with some exponent $\nu\geq 1$.
Then $Y$ is pre-Cartier-split but not Frobenius-split.
\end{proposition}

\proof
Suppose $Y$ is Frobenius split, and choose a retraction for $\O_{Y'}\subset F_*(\O_Y)$.
The latter defines a non-zero global section $s$ of $\omega_{Y/Y'}=\omega_Y^{\otimes 1-p}$.
This sheaf is numerically trivial, so the  map $s:\O_Y\ra\omega_{Y/Y'}$ is bijective.
In turn, the order $p^\nu$ divides $p-1$, contradiction.
Thus $Y$ is not Frobenius split. By assumption we have $c_1=0$, and   
 the dualizing sheaf $\omega_Y$ is not $(p-1)$-torsion in the Picard group.
Consequently, Lemma \ref{not one-dimensional} gives $\Hom(B\Omega^1_Y,\O_{Y'})=0$.
A priori, the connecting map $\Hom(B\Omega^1_Y,\O_{Y'})\ra \Ext^1(\Omega^1_{Y'},\O_{Y'})$ is zero, 
hence the scheme   $Y$ is   pre-Cartier split.
\qed

\medskip
Recall that $n=\dim(Y)$. We now look at the right end of the cochain complex $F_*(\Omega_Y^\bullet)$.
Using $\Omega_{Y'}^n=\omega_{Y'}$ and $Z \Omega^n_Y=F_*(\omega_Y)$ we get a   commutative diagram
\begin{equation}
\label{right four-term sequence}
\begin{tikzcd}[row sep=tiny, column sep=small ]
		&				&0\arrow[dr]			& 		& 0\\		
		&				&						& B\Omega^n_Y\arrow[dr]\arrow[ur]\\
0\arrow[r]	& Z\Omega^{n-1}_Y\arrow[r]	& F_*\Omega^{n-1}_Y\arrow[ru]\arrow[rr, "d"']	& 		& F_*\omega_Y\arrow[r]	& \omega_{Y'}\arrow[r]	& 0,\\
\end{tikzcd}
\end{equation}
where the horizontal four-term sequence is exact,   and obtained by splicing the two short exact sequences with kinks.

\begin{proposition}
\mylabel{surjective map}
Suppose that the order of the dualizing sheaf $\omega_Y$   in the Picard group coincides with the characteristic $p\geq 2$. 
Then the canonical map 
$$
H^{n-1}(Y',\uHom (F_*\Omega_Y^{n-1},\O_{Y'})) \lra  H^{n-1}(Y', \uHom(Z\Omega^{n-1}_Y,\O_{Y'})) 
$$
is surjective.
\end{proposition}

\proof
According to Lemma \ref{dualizing sheaves}, we have $F_*(\omega_Y)=\uHom(F_*\O_Y,\O_{Y'})$, so its dual sheaf
gets identified with $F_*(\O_Y)$, via biduality.
Dualizing the short exact sequence to the right in \eqref{right four-term sequence}, we thus get
$$
0\lra \omega_{Y'}^{\otimes-1}\lra F_*\O_Y\lra \uHom(B\Omega^n_Y,\O_{Y'})\lra 0.
$$
The resulting long exact sequence shows that the induced homomorphism  
$$
H^n(Y,\O_Y)\lra H^n(Y',\uHom(B\Omega^n_Y,\O_{Y'}))
$$
is surjective. The term on the left is zero: 
Serre duality yields $h^n(\O_Y)=h^0(\omega_Y)$, and the latter vanish because $\omega_Y$ is numerically trivial
yet $\omega_Y\not\simeq\O_Y$. 

The short exact sequence to the left in \eqref{right four-term sequence} gives an exact sequence
$$
\Ext^{n-1}(F_*\Omega^{n-1}_Y,\O_{Y'})\lra \Ext^{n-1}(Z \Omega^{n-1}_Y,\O_{Y'})\lra \Ext^n(B \Omega^n_Y,\O_{Y'}).
$$
The term on the right vanishes, as we just saw, hence the mapping on the left is surjective.
\qed

\medskip
For dimension $n=2$, the above yields information for cohomology in degree one. We get the following numerical criterion for surfaces:

\begin{theorem}
\mylabel{non-liftable surface}
Suppose $\dim(Y)=2$, that the order of the dualizing sheaf $\omega_Y$   in the Picard group coincides with
the characteristic $p\geq 2$,
and that $$h^1(\Theta_Y) \geq h^1(\uHom(F_*\Omega^1_Y,\O_{Y'})).$$ Then this inequality is
an equality, and  the scheme $Y$ does not lift to the ring $W_2$.
\end{theorem}

\proof
The scheme $Y$ is pre-Cartier split but not Frobenius-split, according to Proposition \ref{pre-cartier not frobenius}.
With  Proposition \ref{surjective map} we get the estimates
$$
h^1(\Theta_Y)\geq h^1( \uHom(F_*\Omega^1_Y,\O_{Y'}))\geq h^{1}( \uHom(Z\Omega^1_Y,\O_{Y'})).
$$
Thus Theorem \ref{numerical criterion} applies: The above inequalities must be equalities,
and the scheme $Y$ does not lift to the ring $W_2$.
\qed

\medskip
The big advantage of the preceding result is that the differentials in the cochain complex
$F_*\Omega^\bullet_Y$ do not enter anymore.
We merely have to compute the first cohomology of the locally free sheaf $\uHom (F_*\Omega_Y^1,\O_{Y'})$.
The next result tells us that under certain assumptions, this dual of Frobenius pushforward
remains a Frobenius pushforward, which   makes the necessary computations of cohomology feasible.

We say that a quasicoherent sheaf $\shE$ on $Y$ admits an \emph{$F$-descend} if $\shE\simeq F^*(\shE')$ for some quasicoherent sheaf $\shE'$
on $Y'$. By fpqc-descend (\cite{SGA 1}, Expos\'e VIII), this   means that   on the fiber product $Y\times_{Y'}Y$ there is a descend datum
$\varphi:\pr_1^*(\shE)\ra\pr_2^*(\shE)$.
According to  \cite{Katz 1970}, Theorem 5.1 such a descend datum can be interpreted as an  integrable connection $\nabla:\shE\ra\shE\otimes\Omega^1_Y$
with   $p$-curvature zero.

\begin{proposition}
\mylabel{frobenius descend}
Suppose    $\shE$ is a  locally free  sheaf on $Y$  that admits an $F$-descent.
Then  $\underline{\Hom}(F_*(\shE^\vee),\O_{Y'}) = F_*(\shE\otimes\omega_{Y/Y'})$.
If moreover the dualizing sheaf $\omega_X$ is $p$-torsion in the Picard group, we have  
$\underline{\Hom}(F_*(\shE\otimes\omega_Y),\O_{Y'}) = F_*(\shE^\vee )$.
\end{proposition}

\proof
Write $\shE=F^*(\shE')$.  Then also $\shE^\vee=F^*(\shE'^\vee)$, and the 
Projection Formula gives $F_*(\shE^\vee) = F_*(\O_Y)\otimes\shE'^\vee$.
With the relations between tensor products and hom modules (\cite{A 1-3}, Chapter II, \S3,  No.s  1--2), 
together with the definition of relative dualizing sheaves one obtains 
$$
\uHom(F_*(\shE^\vee),\O_{Y'}) = \uHom(F_*\O_Y,\O_{Y'})\otimes\shE'= F_*(\omega_{Y/Y'})\otimes\shE'.
$$
Applying the Projection Formula again, we obtain the first assertion.

Now assume that $\omega_Y$ is $p$-torsion in the Picard group. Then Lemma \ref{dualizing sheaves} and biduality gives 
$$
\uHom(F_*(\omega_Y),\O_{Y'}) = \uHom(\uHom(F_*\O_Y,\O_{Y'}),\O_{Y'}) = F_*(\O_Y),
$$
and we can proceed as in the preceding paragraph.
\qed

\medskip
Note that an invertible sheaf $\shL$ admits $F$-descend if and only if it is $p$-divisible in
the Picard group: If $\shL=\shN^{\otimes p}$, we form the base-change $\shL'=\shN\otimes_kk$ under the
Frobenius map $\lambda\mapsto\lambda^p$, and obtain
$\shL = F_Y^*(\shN) = F^*(\shL')$.
Conversely, if $\shL=F^*(\shL')$, we let $\shN=\shL'\otimes_kk$ under the inverse $\lambda\ra\lambda^{1/p}$ 
of the Frobenius map, and get $\shL=F_Y^*(\shN)=\shN^{\otimes p}$.

%===========================================================
\section{The base of the miniversal deformation}
\mylabel{Base miniversal}

\newcommand{\pure}{{\text{\rm pure}}}
\newcommand{\mixd}{{\text{\rm mixed}}}
We now examine liftability via the miniversal formal deformation, which is also called
the  semi-universal formal deformation, or prorepresentable hull in the terminology of  Schlessinger \cite{Schlessinger 1968}.
Let $Y$ be proper and smooth with $h^0(\O_Y)=1$ over a perfect field $k$ of characteristic $p>0$.
Set 
$$
s=h^1(\Theta_Y)\quadand r=h^2(\Theta_Y),
$$
and let $W=W(k)$ be the ring of Witt vectors.
Let $\foY\ra\Spf(A)$ be the \emph{miniversal formal deformation},
where $A$ is a complete local noetherian 
$W$-algebra with residue field $A/\maxid_A=k$.
Then for every lifting $\foY_B$ of $Y$ over some local Artin $W$-algebra $B$ with residue field $ k$,
there is a homomorphism $A\ra B$ with $\foY_B=\foY\otimes_AB$, and the induced map  of cotangent   spaces
$\maxid_A/(pA+\maxid_A^2)\ra \maxid_B/(pB+\maxid_B^2)$ is unique.
 
\begin{lemma}
\mylabel{base versal deformation}
The base of the miniversal deformation is given by  a ring of the form 
$$
A=W[[T_1,\ldots,T_s]]/(f_1,\ldots,f_r),
$$
where the $f_i$ are formal power series with coefficients from $W$, and their
images in $k[[T_1,\ldots,T_s]]$ have no linear terms.
\end{lemma}

\proof
The argument is parallel to the proof for Proposition 1.5 in \cite{Deligne 1976}. Let me sketch Deligne's arguments for
the convenience of the reader:
Since $Y$ is smooth, the isomorphism classes of deformations over the ring of dual numbers $k[\epsilon]$
correspond to vectors in $H^1(Y,\Theta_Y)$, which has dimension $s$. In turn,   $A$ is a quotient of the formal power series ring $R=W[[T_1,\ldots,T_s]]$
by some ideal $I$ contained in the maximal $\maxid=\maxid_R$.
We have to verify  that the ideal can be generated by $r$  elements. 
Consider the   formal subschemes $S'\subset S''$ inside $S=\Spf(R)$ defined by the ideals
$I\supset \maxid I$. Then the obstruction to extend $\foY\ra S'$ to $S''$ lies in the vector space
$H^2(Y,\Theta_Y\otimes I/\maxid I)= H^2(Y,\Theta_Y)\otimes I/\maxid I$. The tensor factor to the left
has dimension $r$, so we may view the obstruction as an $r$-tuple with entries from $I/\maxid I$.
Choose representatives $f_1,\ldots,f_r\in I$. By construction, $\foY\ra S'$ extends to the formal spectrum of $R/(\maxid I+\sum Rf_i)$.
So by the versal property of $\foY$, the projection $R/(\maxid I+\sum Rf_i)\ra R/I$ admits a retraction, and it
follows that $\maxid I+\sum Rf_i=I$. By the Nakayama Lemma, the elements $f_1,\ldots,f_r\in I$ generate the ideal.
These generators indeed have no linear terms modulo $p$, again because the   
isomorphism classes of deformations over the ring of dual numbers $k[\epsilon]$ corresponds to the $r$-dimensional vector space $H^1(Y,\Theta_Y)$.
\qed

\medskip
Write $\maxid=(p,T_1,\ldots,T_s)$ for the maximal ideal of the formal power series ring $W[[T_1,\ldots,T_s]]$.
Note that $p,T_1,\ldots, T_s$ yield  a basis of the cotangent space $\maxid/\maxid^2$, hence these elements
form a regular system of parameters.
Write the beginning of  the formal power series as 
\begin{equation}
\label{relations versal}
f_i(T_1,\ldots,T_s) \equiv \beta_{i,s+1}p + \sum_{j=1}^s \beta_{ij}T_j \quad \text{modulo $\maxid^2$},
\end{equation}
and let  $\bar{\beta}_{ij}\in k$ be the residue classes of the coefficients $\beta_{ij}\in W$.
This defines two matrices
$$
B_\pure=(\bar{\beta}_{ij})_{\substack{1\leq i\leq r\\ 1\leq j\leq s}}\quadand 
B_\mixd=(\bar{\beta}_{ij})_{\substack{1\leq i\leq r\\ 1\leq j\leq s+1}}
$$
with entries from the field $k$. Checking liftability to $W_2$ now translates into a rank computation:

\begin{lemma}
\mylabel{lifting w2}
The scheme $Y$ lifts to $W_2$ if and only if   $\rank(B_\mixd)=\rank(B_\pure)$.
\end{lemma}

\proof
By the  defining properties of  versal deformations, all   liftings
come from $W$-algebra homomorphism $A\ra W_2$. The latter are given by 
$\mu_j\in pW_2$ satisfying the system of equations $f_i(\mu_1,\ldots,\mu_s)=0$ in the ring $W_2$.
Write these truncated Witt vectors as $\mu_j=(0,\lambda_j)$ with scalars $\lambda_j\in k$.
For any lift $\tilde{\lambda}_j\in W_2$ of $\lambda_j\in W_2/pW_2$  we have $p\tilde{\lambda}_j=\mu_j$.
By abuse of notation, we   may also write $ p\lambda_j$ for this element.
The equations $f_i(p\lambda_1,\ldots, p\lambda_s)=0$ in the ring $W_2$ 
translate into the system of linear equations
$$
B_\pure\cdot\begin{pmatrix}
\lambda_1\\
\vdots\\
\lambda_s
\end{pmatrix}
=
-\begin{pmatrix}
\bar{\beta}_{1,s+1}\\
\vdots\\
\bar{\beta}_{r,s+1}
\end{pmatrix},
$$
over the field $k$, in light of \eqref{relations versal}.
Since  the matrix $B_\mixd $ is obtained from  $B_\pure$ by adjoining the transpose of 
$(\bar{\beta}_{1,s+1},\ldots,\bar{\beta}_{r,s+1})$ as additional column, 
solvability of the above system of linear equations means that the two matrices have the same rank.
\qed

\medskip
We now consider the case where the number of relations is $r=1$. Then we may drop the indices for the formal power series,
and we write $f=f_1$ and $\bar{\beta}_j=\bar{\beta}_{ij}$.

\begin{proposition}
\mylabel{nonlift and regularity}
Let $h^1(\Theta_Y)=s$ be arbitrary, but suppose that  $h^2(\Theta_Y)=1$.
Then the scheme $Y$ does not lift to  the ring $W_2$ if and only if  $\bar{\beta}_1=\ldots=\bar{\beta}_s=0$ and $\bar{\beta}_{s+1}\neq 0$.
In this situation, the complete local ring 
$A=W[[T_1,\ldots,T_s]]/(f)$ is regular of dimension $s=h^1(\Theta_Y)$. 
\end{proposition}

\proof
Our matrix  $ B_\text{\rm mixd}$ becomes the vector $(\bar{\beta}_1,\ldots,\bar{\beta}_s,\bar{\beta}_{s+1})$,
and the first assertion follows from Lemma \ref{lifting w2}.
In the cotangent space $\maxid/\maxid^2$ for the regular local ring $W[[T_1,\ldots,T_s]]$, the classes of 
$p,T_1,\ldots,T_s $ form a basis, and the class of the relation $f\in \maxid$ coincides with the first basis vector.
In turn, the residue class ring $A$ remains regular, with $\dim(A)=(s+1)-1=s$.
\qed

\medskip
We say that $Y$ \emph{formally lifts to characteristic zero} if the canonical map $W\ra A$ is injective.
We then regard this map as an inclusion $W\subset A$. Since $W\cap \Nil(A)=0$, 
there is a minimal prime ideal $\primid\subset A$ so that $W\subset  A/\primid$ remains injective.
In turn, $C=A/\primid$ is a complete local ring  with residue field $k=C/\maxid_C$ that is integral,
and whose field of fractions $\Frac(C)$ has characteristic zero. We thus obtain compatible 
infinitesimal deformations
$\foY_i\ra\Spec(C_i)$ of the scheme $ Y$ over the residue class rings $C_i=C/p^{i+1}C$.

\begin{corollary}
\mylabel{lift and singularity}
Assumptions as in the proposition. Suppose furthermore that 
the scheme $Y$   does not lift to the ring $W_2$.
Then  $Y$ formally lifts to characteristic zero if and only if the complete local ring 
$$
A/pA=k[[T_1,\ldots,T_s]]/(f) 
$$
is singular. In this situation, we have $\dim(A/pA)=s-1$.
\end{corollary}

\proof
The local ring $B=W[[T_1,\ldots,T_s]]$ is factorial. 
According to the proposition, the formal power series $f\in B$ is a prime element,
and we may assume that it is of the form $f=p+g$ with some $g\in \maxid_B^2$.
Regard $p\in B$ as another prime element.
If $(f)=(p)$ the ring $A/pA=A=B/pB$ is regular and the scheme $Y$ does not lift to characteristic zero.
On the other hand, if  $(f)\neq (p)$ then the prime element $p$ remains   a regular element in $A= B/fB$.
In turn, the map $W\ra A$ is injective, so the scheme $Y$ formally lifts to characteristic zero.
By Krull's Principal Ideal Theorem, the ring $A/pA=B/(f,p)$ is of dimension $(s+1)-2=s-1$.
Using the description $A/pA=B/(g,p)$ we see that it has embedding dimension $s$, which means that
$A/pA$ is singular.
\qed

\medskip
Recall that the  proper homomorphic images $V/\maxid_V^i$, $i\geq 1$ of   discrete valuation rings $V$ are
exactly the   local artinian principal ideal rings (\cite{McLean 1973}, Theorem 3.3).
Zariski and Samuel call them \emph{special PIR's} (\cite{Zariski; Samuel 1958}, page 245). 
One could also characterize them as the local noetherian rings of dimension zero and embedding dimension at most one.
We propose to call them \emph{jet rings}.

The situation of  Corollary \ref{lift and singularity} is somewhat paradoxical, and perhaps warrants a brief discussion.
The deformations of the scheme $Y$ are 
\emph{unobstructed  in the following sense}:
For each jet ring quotient   $A/\ideala$   
there is smaller ideal $\ideala'\subsetneqq\ideala$ such that $A/\ideala'$ stays a jet ring quotient.
This is because one finds a regular system of para\-meters $f_1,\ldots,f_s\in A$ with $\ideala=(f_1^l,f_2,\ldots,f_s)$,
and then sets $\ideala'=(f_1^{l+1},f_2,\ldots,f_s)$.
One should view the $\Spec(A/\ideala)\subset\Spec(A/\ideala')$ 
as jets of formal  curves inside the base $\Spf(A)$ of the miniversal deformation.
Note that in order to deform over rings   in which $p\neq 0$, one first has to travel over some infinitesimal
neighborhoods in which $p=0$ holds.
On the other hand, one may regard the deformations of $Y$ as \emph{obstructed}:
For certain discrete valuation rings $V$, some jet ring $V/\maxid_V^i$ is the homomorphic
image of $A$, but $V/\maxid_V^{i+1}$ is not. In fact, one may choose $V=W$ with $i=1$,
or $V=k[[T]]$.

Suppose $R$ is  any local noetherian ring. With respect to our prime $p>0$,
one may define the \emph{absolute ramification index}
$$
e(R) = \sup\{i\in\NN\mid p\cdot 1_A\in \maxid_R^i\}\in\NN\cup\{\infty\}.
$$
By Krull's Intersection Theorem,  $e(R)=\infty$ means that $p\in R$ is the zero element, hence $R$ is an $\FF_p$-algebra.
If   $0<e(R)<\infty$,   the residue field $R/\maxid_R$ has characteristic $p>0$, hence 
all other primes $l\neq p$ become invertible, and we get an extension $\ZZ_{(p)}\subset R$ of local rings.
For integral domains $R$, this means flatness.
For discrete valuations rings $R$, our invariant $e(R)$ is then the usual ramification index.
If complete, the ring $R$ becomes an algebra over the ring $W(k)$ of Witt vectors.
Finally, the condition $e(R)=0$ means that $p\in R$ is invertible as well, which makes $R$ into a $\QQ$-algebra.
We see that the  absolute ramification index $e(A)\geq 0$ for the base of the miniversal formal deformation $\foY\ra\Spf(A)$ yields
an numerical invariant of the scheme $Y=\foY_0$ that reflects liftability.

%===========================================================
\section{Proper flat group schemes}
\mylabel{Group schemes}

 Sometimes, first-order liftings  are already precluded  by   Picard schemes.
The goal of this section is to collect some results in this direction, which
mainly rely on the theory of relative group schemes whose structure morphism is proper.
We start with the following general set-up: Let $R$ be a discrete valuation
ring, with residue field $k=R/\maxid_R$ 
and field of fractions $F=\Frac(R)$.
Let $\foG$ be a relative commutative group scheme whose structure morphism
$\foG\ra\Spec(R)$ is proper and flat, and that the closed fiber $\foG_k=\foG\otimes_Rk$ is
connected.
The Stein factorization gives an  affine scheme $\foH=\Spec\Gamma(\foG,\O_{\foG})$,
which is finite and flat over $R$. Using that global sections commute with flat base-change,
one infers that $\foH$ inherits the structure of relative group scheme,
that the canonical map $h:\foG\ra\foH$ is a homomorphism, and that the closed fiber $\foH_k$
is local.
Note that we do not assume that the structure morphism $f:\foH\ra\Spec(R)$ is cohomologically flat,
such that the equality $\O_{\foH}=h_*(\O_{\foG})$ may not be preserved by base-change.

Write $\foG_F^0=(\foG_F)^0$  for the connected component of the origin for the generic fiber, and 
assume throughout that the reduced parts $\foG_{k,\red}=(\foG_k)_\red$ and 
$\foG_{F,\red}^0=(\foG_F^0)_\red$ are geometrically reduced.
This automatically holds if the residue field $k$ is perfect and the function field $F$
has characteristic zero. The assumption ensures that these reduced parts are subgroup schemes, which 
are connected, smooth and proper, hence abelian varieties.
Write  $\foA_F=\foG_{F,\red}^0$, and let $\foA\subset\foG$ be the
Zariski closure, which is an integral closed subscheme that is proper and flat over $R$.
Using that the formation of closures commutes with flat base-change, and we infer that
$\foA\subset\foG$ is a relative subgroup scheme. Since $\Gamma(\foA_F,\O_{\foA})=F$,
the image of the subgroup scheme  $\foA$ in the group scheme  $\foH$ vanishes. 
We shall see below that
the resulting sequence $0\ra\foA\ra\foG\ra\foH\ra 0$ of proper flat relative group schemes   is ``exact''.
One has to exercise some care to make this precise, because 
the category of   commutative group schemes over $R$ is not abelian.
To do so,  view $\foG$ as a scheme over $\foH$, with respect to the canonical morphism
$h:\foG\ra\foH$. As such, it comes with   an action of the induced relative group scheme
$\foA_{\foH}=\foA\otimes_R\Gamma(\O_{\foG})$.

\begin{proposition}
\mylabel{relative abelian}
Assumptions as above. Then the following holds:
\begin{enumerate}
\item
The relative group scheme $\foA$ is an abelian scheme, and we have 
$\foA_k=\foG_{k,\red}$ as   subgroup schemes of $\foG_k$.
\item 
The projection $h:\foG\ra\foH$ is a torsor for the   action of the induced abelian scheme
$\foA_{\foH}$. 
\item
The resulting sequence of group schemes $0\ra\foG_{k,\red}\ra \foG_k\ra\foH_k\ra 0$ is exact.
\end{enumerate} 
\end{proposition}

\proof
First note that by fpqc descent, we may replace the ground ring $R$ by any extension of  discrete valuation rings.
By passing to the strict localization, it suffices to treat the case that
the residue field $k$ is separably closed and that $R$ is henselian.

We first check with the N\'eron--Ogg--Shafarevich Criterion
that the abelian variety $\foA_F$ has good reduction.
Consider the quotient $\Psi_F=\foG_F/\foA_F$, which is an extension of
the group scheme of components by some local group scheme,
and fix a prime $l>0$ that is relatively prime to the characteristic exponent of the residue field $k=R/\maxid_R$,
and that does not divide the order of $\Psi_F$. Write $O_R\subset\foG$ for the zero section.
The cartesian diagram
$$
\begin{CD}
\foG[l^n]	@>>>	O_R\\
@VVV			@VV\can V\\
\foG		@>>l^n>	\foG
\end{CD}
$$
defines a relative subgroup scheme $\foG[l^n]$. Its formation commutes with base-change in $W$,
and the structure morphism $\foG[l^n]\ra\Spec(R)$ is proper. 
Since multiplication by $l^n$ is finite on 
abelian varieties and finite group schemes, we see that
the structure morphism is quasifinite, hence finite.
Moreover, 
both fibers have length $l^{2ng}$, where $g\geq0 $ is the common dimension of the
two abelian varieties $\foG_{k,\red}$ and $\foA_F$. With \cite{Hartshorne 1977}, Chapter III, Theorem 9.9 
we conclude that the structure morphism $\foG[l^n]$ is locally free of degree $d=l^{2ng}$, that the generic fiber $\foG[l^n]_F$ is contained
in $\foA_F$, and that the closed fiber  $\foG[l^n]_k$ is contained in $\foG_{k,\red}$.
Since the residue field $k$ is separably closed and the ring is henselian, 
we see that the relative group schemes $\foG[l^n]$ are constant. In particular, each $F^\sep$-valued point on $\foA_F[l^n]$ comes from a $F$-valued point.
This ensures that the N\'eron model  $\foA'$  of $\foA_F$ over $R$ is an abelian scheme
(\cite{Bosch; Luetkebohmert; Raynaud 1990}, Chapter 7.4, Theorem 5).

We thus have  mutually inverse birational map $\foA\dashrightarrow\foA'$ and $\foA'\dashrightarrow \foA$, 
whose domains of definitions contain the generic fiber $\foA_F=\foA'_F$.
Choose some integral scheme $\foX$ and some  proper birational  morphisms $\foA\leftarrow \foX\ra\foA'$
over $R$ that become identities over $F$. The union $\bigcup_{n\geq 0} \foA_F[l^n]$ is Zariski dense
in the generic fiber, and the union of the closures become Zariski dense in the closed fibers of $\foA$ and $\foA'$.
It follows that the strict transforms of $\foA_k$ and $\foA'_k$ in $\foX$ coincide.
Consequently, $\foA_k$ is generically reduced. In turn, this closed fiber is
an abelian variety, and $\foA\ra\Spec(R)$ is an abelian scheme. This establishes (i).

To proceed, we employ the theory of stacks.
Consider the relative subgroup scheme $\foA\subset\foG$ and the resulting stack $[\foG/\foA]$.
The latter is  a  fibered category over the category of affine schemes $(\Aff/R)$, with fibers over $U=\Spec(A)$ given
by pairs $(\shT,\varphi)$, where $\shT$ is an $\foA|U$-torsor  and $\varphi:\shT\ra\foG|U$ is
an equivariant morphism. The pairs with $\shT=\foA|U$ and $\varphi$   the  canonical inclusion $\foA|U\subset\foG|U$
define a 1-morphism $\foG\ra [\foG/\foA]$. Note that the  2-fiber product
$\foG\times_{[\foG/\foA]}\foG$ has fiber categories given by  $(g_1,g_2,\psi)$,
where $g_i\in \foG(U)$, and the isomorphism $\psi$ can be viewed as some $a\in\foA(U)$ with $g_2=a+g_1$.
In turn, the canonical  morphism 
\begin{equation}
\label{action projection}
\foA\times\foG\lra \foG\times_{[\foG/\foA]}\foG,\quad (a,g)\longmapsto (ag,g)
\end{equation}
is a 1-isomorphism.
For the abelian scheme $\foA$, the structure morphism $\foA\ra\Spec(R)$ is smooth, separated and   of finite type,
hence $[\foG/\foA]$ is an \emph{Artin stack}, for example by \cite{Laumon; Moret-Bailly 2000}, Example 4.6.1. 

The 1-morphisms  $[\foG/\foA]\ra\Spec(R)$ is separated. To see this,
we apply the Valuative Criterion (loc.\ cit.\, Proposition 7.8):
Suppose $(\shT_1,\varphi_1)$ and $(\shT_2,\varphi_2)$ are two objects of $[\foG/\foA]$ over
the spectrum $U$ of a complete valuation ring $A$ with algebraically closed residue field,
and $\alpha_F:(\shT_1,\varphi_1)_F\ra(\shT_2,\varphi_2)_F$ is an isomorphism over $F=\Frac(A)$.
We have to check that there is at most one extension to an isomorphism over $A$.
The assumptions on $A$ ensure that the torsors admit sections, so  it suffices to treat the case that both $\shT_1,\shT_2$ coincide with $\foA|U$.
Then $\alpha_F$ can be identified with an element of $s_F\in\foA(F)$, and the separatedness of $\foA\ra\Spec(R)$
ensures that there is at most one extension.

The algebraic stack $[\foG/\foA]$  is actually an algebraic space. To see this, we apply loc.\ cit.\ Corollary 8.1.1
and have to verify the following: Suppose $(\shT,\varphi)$ is an object of the stack over some $U=\Spec(A)$,
and $\alpha$ is an automorphism of $(\shT,\varphi)$, then actually $\alpha=\id$.
This problem is local, so it suffices to treat the case $\shT=\foA|U$. Then $\alpha$ is
the translation with respect to some section $s\in \foA(U)$. Since the $\foA$-action on $\foG$ is free,
we must have $s=e$, thus $\alpha=\id$. 
 % Using the Valuative Criterion (loc.\ cit.\, Proposition 7.8),
% one sees that the 1-morphism $[\foG/\foA]\ra\Spec(R)$ is separated.
% With the characterization of algebraic spaces (loc.\ cit.\ Corollary 8.1.1) we infer
% that $[\foG/\foA]$ is 1-isomorphic to an algebraic space $\foX$.

Since \eqref{action projection} is an isomorphism, the projection
$\foG\ra\foX$ is a $\foG_\foX$-torsor. In turn, the structure morphism $\foX\ra\Spec(R)$
is of finite type.
Forming the stack $[\foA/\foG]$ commutes with base-change in $R$,
consequently the fibers of $\foX$ are finite. It follows that $\foX$ is a scheme (loc.\ cit., Theorem A.2).
Since $\foG\ra\Spec(R)$ is proper, $\foG\ra\foX$ is surjective and $\foX\ra\Spec(R)$ is separated and of finite type,
the latter must be proper (\cite{EGA II}, Corollary 5.4.3), hence finite. Since the projection $h:\foG\ra\foX$ is a torsor for the
abelian variety, we see that $\O_{\foX}=h_*(\O_{\foG})$. In turn, $\foG$ and $\foX$ have the same Stein
factorization, and it follows $\foH=\foX$. This establishes (ii).
Since $\foA_k=\foG_{k,\red}$,  and forming the stack $[\foG/\foA]$ commutes with base-change
in $\Spec(R)$, we also have (iii).
\qed

\medskip
We now consider the following more special situation: 
Suppose that $k$ is a perfect field, and $G$ is a commutative group scheme over $k$
that is proper and connected. 
Then the reduced part $G_\red$ is a subgroup scheme that is  an abelian variety,
and the quotient $L=G/G_\red$ is a local group scheme.
Let $L^\mult$ be its multiplicative part, such that $U=L/L^\mult$ is unipotent.
Since $k$ is perfect, the resulting extension splits uniquely, and
we have $L=L^\mult\times U$, see \cite{Demazure; Gabriel 1970}, Chapter IV, \S3, Theorem 1.1.
Note that $L^\mult $ corresponds to finite Galois modules, 
whereas $U$ is given by a Dieudonn\'e module of finite length.
Write $W=W(k)$ for the ring of Witt vectors,   $W_2$ for its truncation of length two,
and $\alpha_{p^n}=\GG_a[F^n]$ for the iterated Frobenius kernel.
 
\begin{proposition}
\mylabel{group schemes no lift}
Suppose the unipotent group scheme $U=L/L^\mult$ contains   $\alpha_{p^n}$ as a direct summand,
for some exponent $n\geq 1$. 
Then   $L$ does not lift to the ring $W_2$, and   $G$ does
not lift to the ring $W$. 
\end{proposition}

\proof
Seeking a contradiction, we assume that there is a relative group scheme $\foG$ whose 
structure morphism $\foG\ra \Spec(W)$ is proper and flat, with closed fiber $\foG_k=G$.
According to Proposition \ref{relative abelian}, there is relative group scheme $\foH$ whose
structure morphism is finite, with closed fiber $\foH_k=L$.
In particular, $L$ lifts to the ring $W_2$, which  reduces the second assertion to the first.

Now suppose that we have a relative group scheme $\foL\ra\Spec(W_2)$ whose structure morphism
is finite, with closed fiber $\foL_k=L$, and consider the ensuing Hopf algebra $H=\Gamma(\foL,\O_{\foL})$.
As a $k$-algebra, the fiber ring  $\bar{H}=H/pH$ takes the form
$\bar{H}=k[T_1,\ldots,T_r]/(T_1^{p^{n_1}},\ldots,T_r^{p^{n_r}})$ for some integer $r\geq 0$ and
some exponents $n_i\geq 0$, 
according to \cite{Demazure; Gabriel 1970}, Chapter III, \S3,  Corollary  6.3. 

By assumption, we have a decomposition $L=\alpha_{p^n}\oplus L'$.
The first factor  is the spectrum of the local Artin ring $k[t]/(t^{p^n})$, with comultiplication
$t\mapsto t\otimes 1+1\otimes t$.
The projection $L\ra\alpha_{p^n}$ corresponds to an inclusion of Hopf algebras $k[t]/(t^{p^n})\subset \bar{H}$,
and we may assume $t=T_1$ and $n=n_1$.

Clearly, the $k$-algebra $\bar{H}$ is a complete intersection,
and the $\bar{H}$-module  $\Omega^1_{\bar{H}/k}$ is freely generated by the  differentials $dT_1,\ldots, dT_r$.
In turn, $\Ext^1(\Omega^1_{\bar{H}},\bar{H})=0$. It follows that all lifts of the scheme $\Spec(\bar{H})$
to the ring $W_2$ are isomorphic, so we may write the $W_2$-algebra as 
$H=W_2[T_1,\ldots,T_r]/(T_1^{p^{n_1}},\ldots, T_r^{p^{n_r}})$.
Using multi-index notation, we observe that
the monomials $T^a= \prod_{i=1}^rT_i^{a_i}$ form a basis for the underlying $W_2$-module, 
with $0\leq a_i<p^{n_i}$.
The comultiplication takes the form
$$
\Delta(t)=t\otimes 1+1\otimes t + p\sum\lambda_{ab}T^a\otimes T^b,
$$
for certain scalars $\lambda_{ab}\in W_2$.
As in   \cite{Oort; Mumford 1968}, first example in the introduction, we now use
that the map $\Delta:H\ra H\otimes H$ is a homomorphism of rings:
On the one hand, thanks to the relation $t^{p^n}=0$ we get $\Delta(t^{p^n})=\Delta(0)=0$.
On the other hand, the  relation $p^2=0$ gives 
$$
\Delta(t)^{p^n}=(t\otimes 1+1\otimes t + p\sum\lambda_{ab}T^a\otimes T^b)^{p^n} = 
(t\otimes 1+1\otimes t )^{p^n}.
$$
With    $t^{p^n}=0$ and the  Binomial Theorem,  this becomes $\sum_{i=1}^{p^n-1}\binom{p^n}{i} t^i\otimes t^{p^n-i}$.
It is well-known that on binomial coefficients of the form $\binom{p^n}{i}$,
the $p$-adic valuation $\nu_p:\ZZ\ra\NN\cup\{\infty\}$
takes the value $n-\nu_p(i)$. In particular, the binomial coefficient  for $i=p^{n-1}$ does not vanish
in the ring $W_2$.
Summing up, we have the  contradiction $0=\Delta(t^{p^n})=\Delta(t)^{p^n}\neq 0$.
\qed
 
\medskip
In light of this, relative  Picard scheme  may preclude  liftings to $W$ or its truncation $W_2$, 
under suitable representability
assumptions. Suppose that $Y$ is a smooth proper scheme over $k$, satisfying $h^0(\O_Y)=1$.   
Let  
$$
b_i=\rank_{\ZZ_\ell}  (\invlim_{\nu} H^i(Y\otimes k^\sep,\mu_{\ell^\nu}^{\otimes i}) ) 
$$
be its $\ell$-adic Betti numbers.

\begin{theorem}
\mylabel{picard no lift}
Set $G=\Pic^0_{Y/k}$. Suppose the local  group scheme $L=G/G_\red$ contains
some $\alpha_{p^n}$, $n\geq 1$ as a direct summand,
and that   $b_1=2(h^1(\O_Y)-h^2(\O_Y))$ holds. Then the scheme $Y$ does not lift to the ring 
$W$.  If moreover $G_\red=0$, the scheme does not even lift to $W_2$.
\end{theorem}

\proof
This relies on some  foundational results on   relative Picard schemes, which we recall first.
Suppose   that $R$ is an arbitrary local noetherian $W$-algebra,
and that $\foY\ra\Spec(R)$ is a proper flat morphism with closed fiber $\foY\otimes_Wk=Y$.
By Artin's result (see \cite{Bosch; Luetkebohmert; Raynaud 1990}, Section 8.3, Theorem 1), 
the relative Picard functor is representable by some relative group space $P=\Pic_{\foY/R}$,
which means a group object in the category of algebraic spaces over $R$.  Moreover, the structure morphism
$P\ra\Spec(R)$ is separated (loc.\ cit.\, Section 8.4, Theorem 3),
and the  condition on the Betti number   ensures that it is also flat
(\cite{Ekedahl; Hyland; Shepherd-Barron 2012}, Proposition 4.2).
The inclusion of $P^\tau=\Pic^\tau_{\foY/R}$ is 
representable by an open and closed embedding (\cite{SGA 6}, Expos\'e XIII, Theorem 4.7 
together with \cite{FGA VI}, Corollary 2.3). Moreover,
the structure morphism $P^\tau\ra\Spec(R)$ is proper
(\cite{Bosch; Luetkebohmert; Raynaud 1990}, Section 8.4, Theorem 4 combined with Theorem 3).

Now suppose that the local ring $R$ is henselian, of dimension $\dim(R)\leq 1$.
According to \cite{Anantharaman 1973}, Theorem 4.B, the algebraic space $P^\tau $ is actually a scheme.
Moreover, the connected components of the scheme $P^\tau$ correspond to 
the connected components of the closed fiber $P^\tau\otimes_Rk$.
Write $\foG\subset P^\tau$ for the connected component with $\foG_k=P^0\otimes_Rk$. This is a subgroup scheme.
By construction, the structure morphism $\foG\ra\Spec(R)$ is proper and flat, and the closed
fiber $G=\foG_k$ is connected.

Seeking a contradiction, we now suppose that $R=W$ is the ring of Witt vectors.
Proposition \ref{relative abelian} applies, and we find some finite flat group scheme $\foH$ over $W$
with closed fiber $\foH_k=G/G_\red$. In particular, the local group scheme $L=G/G_\red$ contains $\alpha_{p^n}$ and
lifts to the ring $W_2$, in contradiction to Proposition \ref{group schemes no lift}.
Finally suppose that $G_\red=0$, and that $R=W_2$ is the ring of truncated Witt vectors.
Then the group scheme $L$ lifts to the ring $W_2$, 
and  Proposition~\ref{group schemes no lift} gives again a contradiction.
\qed

%===========================================================
\section{Vector bundle computations}
\mylabel{Vector bundle}

We now make some computations with vector bundles on surfaces
that will be useful in the following sections.
Suppose our smooth proper scheme $Y$ has dimension $n=2$,
and let $\shE$ be a locally free sheaf. 
The Hirzebruch--Riemann--Roch Theorem $\chi(\shF)=\operatorname{ch}(\shF)\operatorname{td}(\Omega^1_Y)$ 
applied to $\shF=\shE$ and $\shF=\O_Y$ yields the formula
\begin{equation}
\label{hirzebruch riemann roch}
\chi(\shE) = \frac{(D\cdot D)-(D\cdot K_Y)}{2} -c_2(\shE) + \rank(\shE)\chi(\O_Y),
\end{equation}
where  for simplicity we set    $\det(\shE)=\O_Y(D)$ and
$\omega_Y=\O_Y(K_Y)$. Moreover, $c_2(\shE)\in\ZZ$ is the \emph{second Chern number}. Note that this integer is uniquely defined
by the above equation.  

Now suppose that $\shE$ has rank two.
Let $\shE\ra \shF$ be a surjection onto some coherent sheaf that is invertible in codimension one.
Then dual sheaf $\shL=\uHom(\shF,\O_Y)$ is reflexive of rank one (\cite{Hartshorne 1994}, Corollary 1.8),
whence invertible. The canonical map 
$$
\shE=\shE^{\vee\vee} \lra  \shF^{\vee\vee} =\shL^{\vee}
$$
is surjective in codimension one, thus yields an exact sequence $\shE\ra\shL^\vee\ra\O_Z\ra 0$,
where $Z\subset Y$ is a finite    subscheme. 
Let $\shI\subset \O_Y$ be the corresponding coherent ideal sheaf.

\begin{proposition}
\mylabel{description of sequence}
The kernel for the resulting surjection $\shE\ra\shI\shL^\vee$ is isomorphic
to the invertible sheaf $\shL\otimes\det(\shE)$. Moreover, for each point $a\in Z$ the local ring
$\O_{Z,a}$ is of the form $\kappa(a)[[x,y]]/\ideala$ for some parameter ideal $\ideala=(f,g)$.
\end{proposition}

\proof
Let $\shN\subset\shE$ be the kernel in question. The exact sequence
\begin{equation}
\label{long sequence}
0\lra\shN\lra\shE\lra\shL^\vee\lra \O_Z\lra 0
\end{equation}
shows that at each point $a\in Y$, the stalk  $\shN_a$ is a syzygy for
the module $\O_{Z,z}$ over the regular local ring $R=\O_{Y,a}$.
In turn, $\shN$ is locally free.
Taking ranks at the generic point $\eta\in Y$ we see that $\shN$ is invertible.

On the open set $U=Y\smallsetminus Z$, we have $\O_Z|U=0$, thus
$\det(\shE)_U=\shN_U\otimes\shL^\vee_U$. 
This subset $U\subset Y$ contains all points
of codimension one, and with  \cite{Hartshorne 1994}, Theorem 1.12 we deduce that the equality already holds over $Y$.
Thus $\shN=\det(\shE)\otimes\shL$.

Now fix a point $a\in Z$. Set $L=\kappa(a)$, write $\O_{Y,a}^\wedge = L[[x,y]]$ and choose trivializations of $\shE$ and $\shL^\vee$ on some
open neighborhood. The exact sequence \eqref{long sequence} shows that $\O_{Z,a}$ is of the
form $L[[x,y]]/\ideala$ for some ideal $\ideala=(f,g)$. The generators   lie in the maximal ideal and 
form a parameter system, because $\dim(\O_{Z,a})=0$.
\qed

\medskip
We thus have a short exact sequence of coherent sheaves
\begin{equation}
\label{general sequence}
0\lra \shL\otimes\det(\shE)\lra \shE\lra\shI\shL^\vee\lra 0,
\end{equation}
This sequence gives another expression for the second Chern number:

\begin{proposition}
\mylabel{chern number}
In the above situation, we have the formula 
$$
c_2(\shE) + c_1^2(\shL) + c_1(\shL)c_1(\shE) = h^0(\O_Z).
$$
\end{proposition}

\proof
First of all, 
$\chi(\shE) = \chi(\shL\otimes\det(\shE)) + \chi(\shL^\vee) - h^0(\O_Z)$ holds by additivity of Euler characteristics.
Applying Riemann--Roch  this becomes
$$
L^2+(L\cdot D)  + \frac{(D\cdot D) - (D\cdot K_Y)}{2} + 2\chi(\O_Y) - h^0(\O_Z),
$$
where we write  $\shL=\O_Y(L)$ and $\det(\shE)=\O_Y(D)$. 
Together with the Hirzebruch--Riemann--Roch formula \eqref{hirzebruch riemann roch}, this gives the assertion.
\qed
 
\medskip
Under suitable assumptions on the invertible sheaf $\det(\shE)$, $\shL$, $\omega_Y$ 
one obtains formulas for the cohomological invariants of $\shE$:
 
\begin{proposition}
\mylabel{cohomological invariants}
If  in the above situation the dual sheaves for $\shL$ and $\shL\otimes\det(\shE)\otimes\omega_Y^\vee$ have no non-zero global sections,
the cohomological invariants are given by 
$$
h^i(\shE) = \begin{cases}
h^0(\shL\otimes\det(\shE)) 				& \text{for $i=0$;}\\
h^1(\shL\otimes\det(\shE))+h^1(\shL^\vee)+h^0(\O_Z)	& \text{for $i=1$;}\\
h^0(\shL\otimes\omega_Y) 						& \text{for $i=2$.}
\end{cases}
$$
If furthermore $\det(\shE)=\omega_Y$ we obtain the values $h^0(\shE)=h^2(\shE) = h^0(\shL\otimes\omega_Y)$ 
and $h^1(\shE) = 2h^1(\shL\otimes\omega_Y)+ h^0(\O_Z)$.
\end{proposition}

\proof
The groups $H^0(Y,\shI\shL^\vee)\subset H^0(Y,\shL^\vee)$ vanish,  and  the  long exact sequence
for the short exact sequence \eqref{general sequence} gives an identification
$H^0(Y,\shL\otimes\det(\shE))=H^0(Y,\shE)$, which establishes the case $i=0$. 

Likewise get  $H^2(Y,\shE)=H^2(Y,\shI\shL^\vee)$, 
because the group $H^2(Y,\shL\otimes\det(\shE))$ is Serre dual to $H^0(Y,\shL^\vee\otimes\det(\shE)^\vee\otimes\omega_Y)=0$.
The finite subscheme $Z\subset Y$ yields an exact sequence
$0\ra\shI\shL^\vee\ra \shL^\vee\ra\O_Z\ra 0$. In turn, we get a long exact sequence
$$
H^1(Y,\O_Z)\lra H^2(Y,\shI\shL^\vee)\lra H^2(Y,\shL^\vee)\lra H^2(Y,\O_Z).
$$
The outer terms vanish for dimension reason, which establishes the formula for  $i=2$.
For the remaining case, consider the long exact sequence
$$
H^0(Y,\shL^\vee)\lra H^0(Y,\O_Z)\lra H^1(Y,\shI\shL^\vee)\lra H^1(Y,\shL^\vee)\lra H^1(Y,\O_Z).
$$
The outer terms vanish, and we get $h^1(\shI\shL^\vee)=h^0(\O_Z)+h^1(\shL^\vee)$.
By Serre Duality and our assumption on $\shL\otimes\det(\shE)\otimes\omega_Y^\vee$, 
the group $H^2(Y,\shL\otimes\det(\shE))$ vanishes.
From \eqref{general sequence} again we get a short exact sequence
$$
0\lra H^1(\shL\otimes\det(\shE))\lra H^1(Y,\shE)\lra H^1(Y,\shI\shL^\vee)\lra 0,
$$
and the case $i=1$ follows. The formulas for the situation $\det(\shE)=\omega_Y$ are immediate.
\qed

\medskip
A curve $H\subset Y$ is called \emph{ample} or \emph{semiample}  if the invertible sheaf $\O_Y(H)$ has 
the respective property. 

\begin{corollary}
\mylabel{dual invariants}
Suppose that the dualizing sheaf $\omega_Y$ is two-torsion in the Picard group, that $\det(\shE)=\omega_Y$, 
and that $(\shL\cdot H)>0$ for some semiample curve $H\subset Y$. Then  we have
$h^0(\shE^\vee)=h^2(\shE^\vee)=h^0(\shL)$ and $h^1(\shE^\vee) = 2h^1(\shL )+h^0(\O_Z)$.
\end{corollary}

\proof 
We have $\det(\shE^\vee)=\det(\shE)^\vee =\omega_Y^\vee\simeq\omega_Y$.
The wedge product gives a perfect pairing $\shE\otimes\shE\ra\Lambda^2\shE=\omega_Y$,
hence we get identifications $\shE=\uHom(\shE,\omega_Y)$ and  $\shE^\vee=\shE\otimes\omega_Y$.
Tensoring the exact sequence \eqref{general sequence} with $\omega_Y\simeq\omega_Y^\vee$ gives 
$$
0\lra \shN\otimes\omega_Y\lra \shE^\vee\lra \shI\shN^\vee\lra 0
$$
for the invertible sheaf $\shN=\shL\otimes\omega_Y^\vee$. We have $(\shN\cdot H)=(\shL\cdot H)>0$,
so the duals for the invertible sheaves $\shN=\shN\otimes\det(\shE^\vee)\otimes\omega_Y$ have no global
sections.   Proposition \ref{cohomological invariants} applied with $\shE^\vee $ and the above short exact sequence  yields the formulas.
\qed

\begin{corollary}
\mylabel{equality invariants}
Assumptions as in the previous corollary. Then the following are equivalent:
\begin{enumerate}
\item $h^i(\shE)=h^i(\shE^\vee)$ for some degree $0\leq i\leq 2$.
\item $h^i(\shL)=h^i(\shL\otimes\omega_Y)$ for some degree $0\leq i\leq 1$.
\end{enumerate}
If one of these equivalent conditions is true, 
the equalities   hold  for all   $0\leq i\leq 2$.
\end{corollary}

\proof
The previous corollary gives $h^0(\shE^\vee)=h^2(\shE^\vee)= h^0(\shL)$  and $h^1(\shE^\vee)=2h^1(\shL )+h^0(\O_Z)$.
Proposition \ref{cohomological invariants} yields $h^0(\shE)=h^2(\shE)=h^0(\shL\otimes\omega_Y)$ and
$h^1(\shE)=2h^1(\shL\otimes\omega_Y)+h^0(\O_Z)$.
Furthermore, both invertible sheaf $\shL$ and $\shL\otimes\omega_Y$ have no cohomology in degree two, and the
same Euler characteristics.
Suppose we have $h^i(\shL)=h^i(\shL\otimes\omega_Y)$ for some $0\leq i\leq 1$. Then equality holds for all $i\geq 0$,
and so does $h^i(\shE)=h^i(\shE^\vee)$.
Conversely, suppose that we have $h^j(\shE)=h^j(\shE^\vee)$ for some degree $0\leq j\leq 2$.
Then $h^i(\shL)=h^i(\shL\otimes\omega)$ for some $0\leq i\leq 1$, and the assertion follows. 
\qed

%===========================================================
\section{Algebraic surfaces}
\mylabel{Algebraic surfaces}

Let $Y$ be a smooth proper surface over an algebraically closed field $k$
of   characteristic $p=2$.
We now  investigate in what circumstances   Theorem \ref{non-liftable surface} applies, such that the surface $Y$ does not lift
to the ring $W_2$.
Choose some coherent quotient $\Omega^1_Y\ra\shF$ that is invertible in codimension one.
Consider the resulting invertible sheaf $\shL=\uHom(\shF,\O_Y)$ and the ensuing exact sequence
\begin{equation}
\label{filtration sequence omega}
0\lra\shL\otimes\omega_Y\lra \Omega^1_Y\lra \shI\shL^\vee\lra 0,
\end{equation}
where $\shI\subset\O_Y$ is a coherent ideal sheaf corresponding to some finite subscheme $Z\subset Y$.
 
\begin{theorem}
\mylabel{special surface no lift}
Suppose the following three assumptions hold:
\begin{enumerate}
\item The dualizing sheaf $\omega_Y$ has order $p=2$ in the Picard group.
\item There is an semiample curve $H\subset Y$ with $(\shL\cdot H)>0$.
\item The invertible sheaf $\shL\otimes\omega_Y$ is $p$-divisible in the Picard group.
%\item The invertible sheaf $\shL$ or $\shL\otimes\omega_Y$ is $p$-divisible in the Picard group.
\end{enumerate}
Then we have
$$
h^1(\Omega^1_Y)=h^1(\uHom(F_*\Omega^1_Y,\O_{Y'}))\quadand h^0(\shL)\leq h^0(\shL\otimes\omega_Y).
$$
If   the latter inequality is an equality, the  surface $Y$ does not lift to the ring $W_2$.
\end{theorem}

\proof
First, we establish the equality $h^1(\Omega^1_Y)=h^1(\uHom(F_*\Omega^1_Y,\O_{Y'}))$, which is
the main part of the argument.
From  Proposition \ref{cohomological invariants} we know that 
\begin{equation}
\label{numbers}
h^1(\Omega^1_Y)= h^1(\shL\otimes\omega_Y)+h^1(\shL^\vee)+h^0(\O_Z).
\end{equation}
Applying the Frobenius pushforward to \eqref{filtration sequence omega} yields an exact sequence of coherent sheaves
$$
0\lra F_*(\shL\otimes\omega_Y) \lra F_*(\Omega^1_Y) \lra F_*(\shI\shL^\vee)\lra 0,
$$
where the two terms on the left are locally free. In turn, we get an exact sequence of coherent sheaves
\begin{equation}
\label{dual sequence}
0\ra\uHom(F_*(\shI\shL^\vee),\O_{Y'})\ra \uHom(F_*(\Omega^1_Y),\O_{Y'})\ra \uHom(F_*(\shL\otimes\omega_Y),\O_{Y'}) 
\end{equation}
Being duals on a regular two-dimensional scheme, all terms are locally free. The  cokernel for the map on the right is the skyscraper sheaf 
$\shT=\underline{\Ext}^1(F_*(\shI\shL^\vee),\O_{Y'})$.
Moreover, the restriction  map 
$$
\uHom(F_*(\shL^\vee),\O_{Y'}) \lra \uHom(F_*(\shI\shL^\vee),\O_{Y'})
$$
between locally free sheaves is bijective (for example \cite{Hartshorne 1994}, Theorem 1.12). 
% If $\shL$ is $p$-divisible in the Picard group, it admits an   $F$-descent, and
% Proposition \ref{frobenius descend} with $\shE=\shL$ gives 
% $$
% \uHom(F_*(\shL^\vee),\O_{Y'}) = F_*(\shL\otimes\omega_Y)\quadand \uHom(F_*(\shL\otimes\omega_Y),\O_{Y'}) = F_*(\shL^\vee).
% $$
Applying Proposition \ref{frobenius descend} with $\shE=\shL^\vee\otimes\omega_Y$ we get   
$$
\uHom(F_*(\shL\otimes\omega_Y),\O_{Y'}) = F_*(\shL^\vee)\quadand \uHom(F_*(\shL^\vee),\O_{Y'}) = F_*(\shL\otimes\omega_Y).
$$
% Let us treat only the latter case, the former being completely analogous.
The exact sequence \eqref{dual sequence} gives a commutative diagram
\begin{equation}
\label{four-term end sequence}
\begin{gathered}
\begin{tikzcd}[row sep=tiny, column sep=small]
		&					& \phantom{AAA} 0\arrow[dr]					& 		& 0\phantom{}\\		
		&					&								& \shS\arrow[dr]\arrow[ur]\\
0\arrow[r]	& F_*(\shL\otimes\omega_Y)\arrow[r]	& \uHom(F_*(\Omega^1_Y),\O_{Y'})\arrow[ru]\arrow[rr, "d"']	& 	& F_*(\shL^\vee)\arrow[r]	&\shT\arrow[r]	& 0,\\
\end{tikzcd}
\end{gathered}
\end{equation}
where the horizontal four-term sequence is exact, obtained from splicing the two short exact sequences with kinks,
for some coherent sheaf $\shS$.
This gives an inclusion $H^0(Y',\shS) \subset H^0(Y,\shL^\vee)$.
Moreover, we have $H^2(Y',F_*(\shL\otimes\omega_Y))=H^2(Y,\shL\otimes\omega_Y)$, which is dual to $H^0(Y,\shL^\vee)$. 
All these groups vanish, by assumption (ii).
So the long exact sequence for the short exact sequence to the left yields
$$
h^1(\uHom(F_*(\Omega^1_Y),\O_{Y'}))= h^1(\shL\otimes\omega_Y) + h^1(\shS).
$$
On the other hand, the short exact sequence to the right gives 
$$
H^0(Y,\shL^\vee)\ra H^0(Y',\shT)\ra H^1(Y',\shS)\ra H^1(Y,\shL^\vee)\ra H^1(Y',\shT).
$$
The outer terms vanish, by assumption (ii) and for dimension reasons, such that $h^1(\shS)=h^0(\shT)+h^1(\shL^\vee)$.

In light of \eqref{numbers}, it remains to verify $h^0(\shT)=h^0(\O_Z)$. Recall that we started with
an exact sequence $0\ra\shL\otimes\omega_Y\ra\Omega^1_Y\ra\shL^\vee\ra \O_Z\ra 0$, which is a resolution of 
the skyscraper sheaf $\O_Z$ by locally free sheaves. In turn, we get a resolution 
$0\ra F_*(\shL\otimes\omega_Y)\ra F_*(\Omega^1_Y)\ra F_*(\shL^\vee)\ra F_*(\O_Z)\ra 0$ of the skyscraper sheaf
$F_*(\O_Z)$ by locally free sheaves.
Dimension shifting gives
$$
\shT=\underline{\Ext}^1(F_*(\shI\shL^\vee),\O_{Y'}) = \underline{\Ext}^2(F_*(\O_Z),\O_{Y'}).
$$
Since $F:Y\ra Y'$ is finite we have  $h^0(F_*(\O_Z))=h^0(\O_Z)$. Fix a closed point  $b\in Y'$.
It remains to check   that the stalks
\begin{equation}
\label{local duality}
M=F_*(\O_Z)_b\quadand \underline{\Ext}^2(F_*(\O_Z),\O_{Y'})_b=\Ext^2_R(M,R)
\end{equation}
have the same length over the  complete local ring $R=\O_{Y',b}^\wedge$. But this is a general fact:
Let $R/\maxid_R\subset E$ be an injective hull, and $\calC$ be the category of $R$-modules of finite length.
By Matlis Duality (\cite{Matsumura 1989}, Theorem 18.6), the functor  $N\mapsto \Hom_R(N,E)$ induces an anti-equivalence of $\calC$, 
in particular $N$ and $\Hom_R(N,E)$ have the same length. 
Local Duality gives $\Hom_R(\Ext^2_R(N,R),E)=H^0_\maxid(N)=N$, 
with the  2-dimensional local Gorenstein ring $R$ and the finite $R$-module $N$ (see for example \cite{Hartshorne 1967}, Theorem 6.3).
Summing up, the modules in \eqref{local duality} have the same length, and therefore
 $h^1(\Omega^1_Y)=h^1(\uHom(F_*\Omega^1_Y,\O_{Y'}))$.
 
Next, we establish the inequality $h^0(\shL)\leq h^0(\shL\otimes\omega_Y)$.
By \eqref{numbers} and Corollary \ref{dual invariants} we have
\begin{equation}
\label{theta and omega}
h^1(\Omega^1_Y) = 2h^1(\shL^\vee) + h^0(\O_Z)\quadand
h^1(\Theta_Y) = 2h^1(\shL) + h^0(\O_Z).
\end{equation}
Assumption (i) gives
$\chi(\shL)=\chi(\shL\otimes\omega_Y)$, whereas assumption (ii) ensures that
$h^2(\shL)=h^2(\shL\otimes\omega_Y)=0$.
Seeking a contradiction, we now assume $h^0(\shL)> h^0(\shL\otimes\omega_Y)$.
Then we also have $h^1(\shL)>h^1(\shL\otimes\omega_Y)=h^1(\shL^\vee)$,
and with the equations in  \eqref{theta and omega} we obtain $h^1(\Theta_Y^1)> h^1(\Omega^1_Y) =  h^1(\uHom(F_*\Omega^1_Y,\O_{Y'}))$. 
But this contradicts Theorem \ref{non-liftable surface}. Note that for this step we need the assumption that the dualizing sheaf
has order two.

Finally, suppose we have $h^0(\shL)= h^0(\shL\otimes\omega_Y)$. From  Corollary \ref{equality invariants}
we get  $h^1(\Theta_Y^1) = h^1(\Omega^1_Y)=h^1(\uHom(F_*\Omega^1_Y,\O_{Y'}))$, 
and Theorem \ref{non-liftable surface} tells us that the scheme
$Y$ does not lift to the ring $W_2$.
\qed

\medskip
Now suppose we have a \emph{quasielliptic fibration} $f:Y\ra B$. This means that  $B$ is a smooth proper curve,
and the generic fiber $Y_\eta$ is a twisted form of the \emph{rational cuspidal curve}
$\Spec k[t^2,t^3]\cup\Spec k[t^{-1}]$. The fibration gives a short exact sequence 
$$
0\lra f^*(\Omega^1_Y)\lra \Omega^1_Y\lra \Omega^1_{Y/B}\lra 0.
$$
The map on the left is indeed injective, because the function field extension $k(B)\subset k(Y)$ is separable. It follows that
the coherent sheaf $\Omega^1_{Y/B}$ has rank one.
Write $\Omega^1_{Y/B}\ra \shF$ for the quotient modulo the torsion subsheaf.
Then $\shF$ is invertible in codimension one, and we obtain a short exact sequence
$$
0\lra\shL\otimes\omega_Y\lra \Omega_Y^1\lra \shI\shL^\vee\lra 0
$$
attached to the quasielliptic fibration. In order to apply Theorem \ref{special surface no lift},
one has to check that the invertible sheaves  $\shL$ and $\omega_Y$ have certain properties.
Write $F=k(B)$ for the function field of the curve,
and $\Sing(Y_F/F)$ be the scheme of non-smoothness, as defined in \cite{Fanelli; Schroeer 2020}, Section 2.
This is the   closed subscheme of the generic fiber   defined by the 
first Fitting ideal of $\Omega^1_{Y_F/F}$.

\begin{proposition}
\mylabel{fitting ideal}
As Cartier divisors on the generic fiber, we have $\Sing(Y_F/F)=2\xi$
for some closed point $\xi$, and the field extension $F\subset\kappa(\xi)$
is purely inseparable of degree $p=2$.  
\end{proposition}

\proof
We first make an explicit computation with the cuspidal rational curve over $F$. 
The coordinate ring of the first chart is isomorphic to $F[x,y]/(y^2-x^3)$, by setting $x=t^2$ and $y=t^3$.
The  module of K\"ahler differentials is generated by
$dx$ and $dy$ modulo  $x^2dx$. Hence $\Sing(C/F)$ is defined by an additional relation $x^2=0$.
It becomes the spectrum of $F[x,y]/(y^2,x^2)$,
which is radical of length four. In turn, $\Sing(Y_F/F)$ is radical of length four.
It contains no rational point by \cite{Fanelli; Schroeer 2020}, Corollary 2.6. 
Since the field $F$ has $p$-degree $\operatorname{pdeg}(F)=1$, the scheme of non-smoothness has
 residue field $\kappa(\xi)=F^{1/p}$, which has degree two.
Our assertion follows.
\qed

\medskip
The closure $C=\overline{\{\xi\}}$  inside  the quasielliptic surface $Y$ is called the \emph{curve of cusps}.

\begin{proposition}
\mylabel{all fibers simple}
Suppose that all closed fibers $f^{-1}(b)$ are simple, with Kodaira symbol $\II$.
Then we have $\shL=\O_Y(2C)\otimes f^*(\shN)$ for some invertible sheaf $\shN$ on $B$.
\end{proposition}

\proof
By assumption, all geometric fibers in question  are   rational cuspidal curves
$\Spec k[t^2,t^3]\cup\Spec k[t^{-1}]$.
The sheaf of K\"ahler differentials modulo torsion is invertible, and generated on
the first chart by $dt^3$, and on the second chart by $dt^{-1}$.
On the overlap we have $dt^{-1}=t^{-2}dt=t^{-4}dt^3$, which gives the cocycle $t^{-4}\in k[t^{\pm 1}]^\times$.
Its inverse is given by $t^4$, and the resulting divisor coincides with the locus of non-smoothness.

Consider the invertible sheaf $\shM=\shL(-2C)$. The restrictions to   fibers $f^{-1}(b)$, $b\in B$,
are trivial, by the above computation. The direct image $f_*(\O_Y)$ commutes with arbitrary base-change.
By the Theorem of Formal Functions, the direct image $\shN=f_*(\shM)$ is invertible.
According to  the Projection Formula, the adjunction map $f^*(\shN)\ra \shM$ is bijective.
\qed

\medskip
We record the following immediate consequence:

\begin{corollary}
\mylabel{divisibility}
Assumptions as in the proposition. If furthermore  the selfintersection numbers
$(\shL\cdot\shL)$ and $C^2$ vanish,
then $\shL$ is $p$-divisible in the Picard group.
\end{corollary}

\proof
According to the proposition, the invertible sheaf $\shL$ comes from a divisor of the
form $2C+F$, where $F=\sum m_if^{-1}(b_i)$ is a linear combination of fibers.
From $(\shL\cdot \shL) = 4C^2+4C\cdot F$ we infer $C\cdot F=0$. In turn, the divisor
$\sum m_ib_i$ on the curve $B$ has degree zero.
But the group of rational points on the abelian variety  $\Pic^0_C$ is   $n$-divisible for 
any integer $n\geq 1$. It follows that $\shL=\O_Y(2C+F)$ is two-divisible in $\Pic(Y)$.
\qed

%===========================================================
\section{Enriques surfaces}
\mylabel{Enriques surfaces}

Let $k$ be an algebraically closed ground field.
Recall that a smooth surface $Y$ with $h^0(\O_Y)=1$ is called an
\emph{Enriques surface}  if $c_1=0$ and $b_2=10$.
We refer to the monograph of Cossec and Dolgachev \cite{Cossec; Dolgachev 1989} for a comprehensive account.
The group scheme $P=\Pic_{Y/k}^\tau$ of numerically trivial invertible sheaves
is finite of order two, and its group of rational points is generated by the canonical class $K_Y$.
The canonical covering $\epsilon:X\ra Y$ is a torsor under the Cartier dual $G=\uHom(P,\GG_m)$,
and its total space is integral, with cohomological invariants $h^1(\O_Y)=0$ and $h^2(\O_Y)=1$, with  $\omega_Y=\O_Y$.
In characteristic $p\geq 3$, the canonical covering  is a smooth K3 surface, and the base
of the miniversal formal deformation $\foY\ra\Spec(A)$ is given by the ring 
$A=W[[T_1,\ldots,T_{10}]]$.

From now on, we suppose  the characteristic is $p=2$.
Then there are three possibilities for the group scheme $P=\Pic_{Y/k}^\tau$, namely
$\mu_2$ or  $\ZZ/2\ZZ$ or $\alpha_2$. The respective Enriques surfaces $Y$ are aptly called 
\emph{ordinary}, \emph{classical} and \emph{supersingular}.
Ordinary Enriques surfaces behave as in odd characteristics.
For classical and supersingular Enriques surfaces, the group scheme $P$ is unipotent, its
Cartier dual  $G$ is local,  the canonical covering $X$ is singular, and both schemes have trivial
fundamental group. We then say that $Y$ is a \emph{simply-connected Enriques surface},
and $X$ is called the \emph{K3-like covering}.

Let $Y$ be a simply-connected Enriques surface, and $X'\ra X$ be the normalization of 
the K3-like covering. Ekedahl and Shepherd-Barron \cite{Ekedahl; Shepherd-Barron 2004} showed that the ramification divisor 
for the normalization is the preimage of a curve $C\subset Y$ called the \emph{conductrix}. They call
$Y$ an \emph{exceptional Enriques surface} if the \emph{biconductrix} $2C\subset Y$ has  $h^1(\O_{2C})\neq 0$, and give
a beautiful classification of these surfaces in terms of the multiplicities $m_i\geq 1$ and
intersection matrix $(C_i\cdot C_j)$ for the conductrix $C=\sum m_iC_i$,
and also by properties of the Hodge ring $\bigoplus H^i(Y,\Omega_Y^j)$.
Exceptional Enriques surfaces are a priori simply-connected, and both classical and supersingular cases
do occur. The following fact shows that exceptional and supersingular Enriques surfaces
share an important property:

\begin{proposition}
\mylabel{tangent invariants}
The cohomological invariants for the tangent sheaf   of an Enriques surface $Y$ are given by
the following table:
$$
\begin{array}{llll}
\toprule
					& h^0(\Theta_Y)	& h^1(\Theta_Y)	& h^2(\Theta_Y)\\
\midrule
\text{\rm exceptional/supersingular}	& 1		& 12		& 1\\

\text{\rm otherwise}			& 0		& 10		& 0\\
\bottomrule
\end{array}
$$
\end{proposition}

\proof
First note that $\chi(\Theta_Y)= -c_2 + 2\chi(\O_Y) = 10$.
Furthermore, we have $\Theta_Y=\Omega_Y^1\otimes\omega_Y$, and Serre Duality gives
$h^i(\Theta_Y)=h^{2-i}(\Theta_Y)$.
So it suffices to verify the values in degree $i=0$.
Ordinary Enriques surfaces have $h^0(\Theta_Y)=0$, whereas supersingular have $h^0(\Theta_Y)=1$,
according to \cite{Cossec; Dolgachev 1989}, Proposition 1.4.2.
Now suppose that $Y$ is classical. By \cite{Ekedahl; Shepherd-Barron 2004} 
the Enriques surface $Y$ admits non-zero global vector fields if and only if $Y$ is exceptional, and then
$h^0(\Theta_Y)=1$.
\qed
 
\medskip 
We now come to the main result of this paper.  

\begin{theorem}
\mylabel{versal deformation enriques}
Let $Y$ be an Enriques surface, and $\foY\ra\Spf(A)$ be its  miniversal formal deformation.
Then the  complete local noetherian ring  $A$   is regular with   $11\leq \dim(A)\leq 12$,
and flat as $W$-algebra.
Moreover, the following are equivalent:
\begin{enumerate}
\item The Enriques surface $Y$ is exceptional or supersingular.
\item The scheme $Y$ does not lift  to the ring $W_2$.
\item The absolute ramification index is $e(A)\geq 2$.
\item The dimension is $\dim(A)=12$.
\end{enumerate}
\end{theorem}

\proof
Suppose first that $Y$  is neither exceptional nor supersingular. Then we have $h^1(\Theta_Y)=10$ and $h^2(\Theta_Y)=0$,
hence $A=W[[T_1,\ldots,T_{10}]]$, and the assertion is immediate.

Now suppose that $Y$ is exceptional or supersingular, such that $h^1(\Theta_Y)=12$ and $h^2(\Theta_Y)=1$.
According to Proposition \ref{nonlift and regularity}, we merely have to check that the scheme $Y$ does not lift to the ring $W_2$.
For supersingular Enriques surfaces, this follows from \cite{Ekedahl; Hyland; Shepherd-Barron 2012}, Proposition 4.6,
compare also Theorem \ref{picard no lift}, and also \cite{Shepherd-Barron 2017}, Theorem 7.1 for a description of the base
of the versal deformation in terms of invariant rings.

It remains to treat the case that $Y$ is classical and exceptional.
To show that $Y$ does not lift to the ring $W_2$ 
we now  check that the assumptions of Theorem \ref{special surface no lift} are satisfied.
Since $Y$ is classical, the dualizing sheaf $\omega_Y$ has order $p=2$ in the Picard group.
Let $C\subset Y$ be the conductrix, and consider the invertible sheaf $\shL=\omega_Y(2C)$, such that
$\shL\otimes\omega_Y=\O_Y(2C)$. Obviously, $\shL\otimes\omega_Y$ is $p$-divisible in the Picard group,
and $(\shL\cdot H)>0$ for every ample curve $H\subset Y$.
According to the proof of Proposition 0.5 in \cite{Ekedahl; Shepherd-Barron 2004}, there is a short exact sequence
$$
0\lra \shL\otimes\omega_Y\lra \Omega^1_Y\lra \shI\shL^\vee\lra 0,
$$
where  $\shI$ is the ideal sheaf of some   finite subscheme $Z\subset X$.
We have $h^0(\Omega^1_Y)=1$, and the above short exact sequence
immediately gives  $h^0(\shL\otimes\omega_Y)=1$.
Tensoring with $\omega_Y$ we obtain another the exact sequence
$0\ra \shL\ra \Theta_Y\ra \shI\shL^\vee\otimes\omega_Y\ra 0$, which gives an
exact sequence
$$
0\lra H^0(Y,\shL)\lra H^0(Y,\Theta_Y)\lra H^0(Y,\shI\shL^\vee\otimes\omega_Y).
$$
Clearly, the term on the right vanishes. Moreover, we have $h^0(\Theta_Y)=1$ by Proposition \ref{tangent invariants}
and conclude $h^0(\shL)=1$. In particular the equality $h^0(\shL)=h^0(\shL\otimes\omega_Y)$ holds.
We thus may apply Theorem \ref{special surface no lift} and get that the scheme $Y$ does not
lift the ring $W_2$.
\qed

\medskip
Note that the Hodge--de Rham spectral sequence $E_1^{rs}=H^s(Y,\Omega_Y^r)\Rightarrow H^{r+s}(Y,\Omega_Y^\bullet)$
degenerates on the $E_1$-page if and only if the Enriques surface $Y$ is not supersingular,
according to \cite{Illusie 1979}, Proposition 7.3.8. In particular, \cite{Deligne; Illusie 1987}, Corollary 2.4 does not apply for
such Enriques surfaces. Moreover, I do not see any ample invertible sheaf   violating
the Kodaira--Akizuki--Nakano Vanishing $H^r(Y,\Omega_Y^s\otimes\shL)=0$ for $r+s>2$, so loc.\ cit. Corollary 2.8
also does not help to establish non-existence of liftings to $W_2$.

%===========================================================
\section{Bielliptic surfaces}
\mylabel{Bielliptic surfaces}

Let $k$ be an algebraically closed field, and $Y$ be a smooth proper $k$-scheme with $\dim(Y)=2$ and  $h^0(\O_Y)=1$. 
Let us say that $Y$ is a \emph{bielliptic surface} if $c_1=0$ and $b_2=2$.
By the   Enriques classification according to Bombieri and Mumford (\cite{Bombieri; Mumford 1977} and \cite{Bombieri; Mumford 1976}),
we then have 
$$
b_1=b_3=2\quadand c_2= \chi(\O_Y)=0.
$$
Moreover, the  number  $h^1(\O_Y)=h^2(\O_Y)+1$ is either one or two. Note that Bombieri and Mumford  
used the terms \emph{hyperelliptic} and \emph{quasi-hyperelliptic surfaces}.

Throughout,   $Y$ denotes a bielliptic surfaces. Then
$Y=(E\times C)/G$,  
where the first factor $E$ is  elliptic curve, the second factor $C$ is either another elliptic curve or the rational cuspidal curve
$\Spec k[t^2,t^3]\cup\Spec k[t^{-1}]$,
and the finite group scheme  $G$ acts diagonally via   inclusions $G\subset E$ and $G\subset\Aut_{C/k}$.
The action is free on $E$ but non-free on $C$,    and the possible orders   $\ord(G)=h^0(\O_G)$ are the numbers  $d=2,3,4,6$.
The two projections for the product $X=E\times C $ induce two fibrations
$$
B=E/G\stackrel{f}{\longleftarrow} Y  \stackrel{g}{\lra}   C/G=\PP^1
$$
on the quotient $Y=X/G$, where $B$ is another elliptic curve.
Both projections are genus-one fibrations, and $f:Y\ra B$ is quasielliptic if and only
if $C$ is the rational cuspidal curve. In this case, we are in characteristic $p=2$ or $p=3$, and 
the   group scheme $G$ is non-reduced.
Moreover, all closed fibers $f^{-1}(b)$ are simple with Kodaira symbol $\II$.

Consider the one-dimensional representation $\chi:G\ra\GL(H^0(C,\omega_C))=\GG_m$,
and the invertible sheaf 
$\shN=\uHom(R^1f_*(\O_Y),\O_B)$ on the elliptic curve $B=E/G$.

\begin{proposition}
\mylabel{order dualizing}
The dualizing sheaf is $\omega_Y=f^*(\shN)$. Moreover, the common order of $\omega_Y$ and $\shN$ in the
Picard groups coincides with  the order for the subgroup scheme $\chi(G)\subset \GG_m$.
\end{proposition}

\proof
The Canonical Bundle Formula (\cite{Bombieri; Mumford 1977}, Theorem 2) gives $\omega_Y=f^*(\shN)$.
By the Projection Formula, the sheaves $\omega_Y$ and $\shN$ have the same order.
If $C$ is an elliptic curve, the assertion on the order of $\omega_Y$ is given
in loc.\ cit., page 37.
If $C$ is the rational cuspidal curve  the proof for Proposition 8 in \cite{Bombieri;  Mumford 1976} gives the assertion.
\qed

\medskip
Let us now examine the   quasielliptic situation in the most important case  $p=\ord(G)=2$ in more detail:

\begin{proposition}
\mylabel{bielliptic invariants}
Suppose that   $p=2$, and that the bielliptic surface $Y=(E\times C)/G$ is formed  with
the rational cuspidal curve $C$ and the group scheme $G=\mu_2$. Then
$h^1(\O_Y)=1$ and $h^2(\O_Y)=0$, and the cohomological invariants
for the tangent and cotangent sheaves are 
$$
h^i(\Theta_Y) = \begin{cases}
1	& \text{for $i=0,2$};\\
2	& \text{for $i=1$}
\end{cases}
\quadand
h^i(\Omega^1_Y) = \begin{cases}
3	& \text{for $i=0,2$};\\
6	& \text{for $i=1$.}
\end{cases}
$$
\end{proposition}

\proof
Since $G=\mu_2$ is simple and acts non-trivially on $C$, the representation $\chi:G\ra\GL(H^0(C,\O_C))$
is a monomorphism, so the  dualizing sheaf  has order two. It follows that $h^2(\O_Y)=0$ and hence  $h^1(\O_Y)=1$.
Recall that $f:Y\ra B$ denotes the quasielliptic fibration.
Since $\omega_Y=f^*(\shN)$, we see that $\shN$ has order two in $\Pic(B)$.
The inclusion of group schemes $G\subset E$ shows that the elliptic curve $E$ is ordinary, and the same
holds for the isogeneous curve $B=E/G$. 
Note that up to isomorphism, $\shN$ is the only invertible sheaf of order two.

The faithful  action of the height-one group scheme  $G=\mu_2$ on the rational cuspidal curve $C=\Spec k[t^2,t^3]\cup\Spec k[t^{-1}]$
corresponds to a non-zero vector field $\delta\in H^0(C,\Theta_C)$ satisfying $\delta^{[2]}=\delta$.
Write $D_{t^{-1}}$ for the derivative with respect to the variable $t^{-1}$.
As explained in \cite{Schroeer 2007}, Section 3 the Lie algebra $\mathfrak{g}=H^0(C,\Theta_C)$ is four-dimensional, 
and we can write  $\delta=P(t^{-1})D_{t^{-1}}$ for some   polynomial of the form
$$
P(t^{-1}) = \lambda_4t^{-4} + \lambda_2t^{-2} +\lambda_0 + \lambda_1t^{-1}
$$
The condition $\delta^{[2]}=\delta $ means $\lambda_1\neq 0$, and the condition 
that the singularity of $C$ is not a fixed point means $\lambda_4\neq 0$.
The polynomial   is separable, because its derivative is $P'(t^{-1})=\lambda_1$. 
As explained in \cite{Schroeer 2007}, Section 1   (compare also \cite{Kondo; Schroeer 2019}, Section 3)
its four roots define the fixed scheme $C^G$.
Let $c_1,\ldots,c_4\in \PP^1=C/G$ be the images of the fixed points,
and write $0\in\PP^1$ for the image of the singularity $0\in C$.
It follows that the   $ g^{-1}(c_i)$ are precisely the multiple fibers.
These fibers are tame, with multiplicity   $m=2$. Write $B_i=g^{-1}(c_i)_\red$ for the corresponding half-fibers,
and also  set $E_0=g^{-1}(0)$.
We have chosen this notation because the canonical morphisms $Y\ra B=E/G$ and $E\times C\ra Y$ induce
identifications $B_i=B$ and $E\times\{0\}=E_0$. 
From the Canonical Bundle Formula we get $\omega_Y=\O_Y(-2E_0 + B_1+\ldots+B_4)$.

Since all fibers for the quasielliptic fibration $f:Y\ra B$ are simple with Kodaira symbol $\II$,
the coherent sheaf $\Omega^1_{Y/B}$ modulo torsion  is invertible. Setting $\shL=\uHom(\Omega^1_{Y/B},\O_Y)$
we obtain a short exact sequence
\begin{equation}
\label{sequence bielliptic}
0\lra\shL\otimes\omega_Y\lra\Omega_Y^1\lra\shL^\vee\lra 0.
\end{equation} 
Note  that $E_0=g^{-1}(0)$ is the curve of cusps.
According to Proposition \ref{all fibers simple}, we have $\shL=\O_Y(2E_0)\otimes f^*(\shN')$
for some invertible sheaf $\shN'$ on the elliptic curve $B=E/G$. 
We have $E_0^2=0$, and with  $c_1=c_2=0$ and Proposition \ref{chern number} we also get $(\shL\cdot\shL)=0$.
It follows that $\shN'$ is numerically trivial.
We claim that it is has order two, such that $\shN'=\shN$ and $\shL\otimes\omega_Y=\O(2E_0)$.
Consider the half-fiber $B_1=g^{-1}(c_1)_\red$ and the resulting exact sequence
$$
0\lra \O_{B_1}(-B_1)\lra \Omega^1_Y|B_1\lra \Omega_{B_1}^1\lra 0.
$$
With respect to the identification $B_1=B$, the outer terms are $\shN$ and $\O_B$.
Since $\Ext^1(\O_B,\shN)=0$, we obtain $\Omega^1_Y|B_1=\shN\oplus\O_B$.
Restricting the  short exact sequence \eqref{sequence bielliptic} to the curve $B_1$ gives a surjection
$\shN\oplus\O_B=\Omega^1_Y|B_1\ra \shN'^\vee $. Since both $\shN$ and $\shN'$ have degree zero, it follows
that either $\shN'=\shN$ or $\shN'=\O_B$.

Seeking a contradiction, we suppose $\shN'=\O_B$.
Then $\shL=\O_Y(2E_0)$, and the projection formula for the elliptic fibration
$g:Y\ra \PP^1$ gives
$h^0(\shL) = h^0(\O_{\PP^1}(2))=3$.
Furthermore, we have $\shL\otimes\omega_Y=\O_Y(B_1+\ldots+B_4)$.
Each global section vanishes only along curves that are vertical
with respect to $g:Y\ra \PP^1$. Using $\O_{B_i}(B_i)\neq \O_{B_i}$,
we infer $h^0(\shL\otimes\omega_Y)=1$.
Now $h^0(\shL)=3>1=h^0(\shL\otimes\omega_Y)$ contradicts Theorem~\ref{special surface no lift}.

In turn, we have $\shN'=\shN$, hence $\shL\otimes\omega_Y=\O_Y(2E_0)$ and
$\shL=\O_Y(B_1+\ldots+B_4)$. As above, this gives $h^0(\shL)=1$ and $h^0(\shL\otimes\omega_Y)=3$. 
Now
$$
h^0(\Omega_Y^1)=h^2(\Omega_Y^1)=h^0(\shL\otimes\omega_Y)=3\quadand
h^0(\Theta_Y^1)=h^2(\Theta_Y^1)=h^0(\shL)=1
$$
follow  from Proposition \ref{cohomological invariants} and Corollary \ref{equality invariants}. 
The Hirzebruch--Riemann--Roch Formula gives $\chi(\Omega^1_Y)=\chi(\Theta_Y)=0$,
and the values in degree $i=1$  follow as well.
\qed

\medskip
Note that we cannot deduce  non-liftability from  Theorem \ref{special surface no lift}, because
the inequality $h^0(\shL)\leq h^0(\shL\otimes\omega_Y)$ is not an equality.
The situation changes if the  group scheme $G$ is unipotent.
To simplify notation, we now write $\O_Y(n)=g^*(\O_{\PP^1}(n))$ for the pullback of invertible
sheaves along the elliptic fibration $g:Y\ra \PP^1=C/G$.

\begin{proposition}
Suppose that   $p=2$, and that the bielliptic surface $Y=(E\times C)/G$ is formed  with
the rational cuspidal curve $C$ and the group scheme $G=\alpha_2$.
Then the cotangent sheaf sits is a non-split short exact sequence
\begin{equation}
\label{bielliptic extension}
0\lra\O_Y(2)\lra\Omega^1_Y\lra\O_Y(-2)\lra 0,
\end{equation}
in particular we have $\omega_Y=\O_Y$ and $\Theta_Y=\Omega^1_Y$. The
cohomological invariants are given by the formulas  $h^1(\O_Y)=2$, $h^2(\O_Y)=1$ and
$$
h^i(\Omega^1_Y)=h^i(\Theta_Y)  = \begin{cases}
3	& \text{for $i=0,2$};\\
6	& \text{for $i=1$.}
\end{cases}
$$
\end{proposition}
 
\proof
Now the representation $\chi:G\ra\GL(H^0(C,\omega_C))$ is trivial,
whence $\omega_Y=\O_Y$  by Proposition \ref{order dualizing}, such that   $h^2(\O_Y)=1$ 
and  $h^1(\O_Y)=2$. Moreover, $\Theta_Y=\Omega_Y^1\otimes\omega_Y$
is isomorphic to $\Omega^1_Y$. Again set $\shL=\uHom(\Omega^1_{Y/B},\O_Y)$,
and consider the short exact sequence 
\begin{equation}
\label{sequence unipotent}
0\lra \shL\lra\Omega^1_Y\ra \shL^\vee\lra 0
\end{equation}
stemming from the quasielliptic fibration $f:Y\ra B$.
As in the preceding proof, we have $\shL=\O_Y(2E_0)\otimes f^*(\shN')$ for some
numerically trivial sheaf $\shN'$ on the elliptic curve $B=E/G$. We claim that in the present situation 
$\shN'\simeq\O_B$. Choose a fixed point  or the $G$-action on $C$, let $b_1\in \PP^1$ be its image,
and  $B_1=g^{-1}(b_1)_\red$ be the resulting copy of $B$. The Adjunction Formula
shows that  the conormal sheaf $\O_{B_1}(B_1)$ is trivial.
In turn, both outer terms in the short exact sequence
$0\ra \O_{B_1}(-B_1)\ra\Omega^1_Y|B_1\ra \Omega^1_{B_1}\ra 0$
are isomorphic to $\O_{B_1}$. Restricting \eqref{sequence unipotent} to $B_1$ gives an inclusion 
$\shN'\subset \Omega^1_Y|B_1$. We infer $\Hom(\shN',\O_B)\neq 0$, and hence $\shN'\simeq \O_B$.
This gives the short exact sequence \eqref{bielliptic extension}.
The long exact cohomology sequence and the Projection Formula yields $h^0(\Omega^1_Y)=h^0(\O_Y(2))= 3$.
By Serre duality   $h^2(\Omega^1_Y)=h^0(\Theta_Y\otimes\omega_Y) = h^0(\Omega^1_Y)=3$.
The Hirzebruch--Riemann--Roch Formula ensures $\chi(\Omega^1_Y)= 0$, and thus $h^1(\Omega^1_Y)=6$.

It remains to check that the extension \eqref{bielliptic extension} does not split.
For this we first   compute   $R^1g_*(\O_Y)$, which can be written
as $R^1g_*(\O_Y)=\O_{\PP^1}(d)\oplus\shF$ for some integer $d$ and some finite sheaf $\shF$.
As explained in the previous proof, the action of $G=\alpha_2$ on the rational cuspidal curve $C=\Spec k[t^2,t^3]\cup\Spec k[t^{-1}]$
is given by some derivation $\delta=P(t^{-1})D_{t^-1}$
with $P(t^{-1})= \lambda_4t^{-4} + \lambda_2t^{-2} +\lambda_0$ with $\lambda_4\neq 0$.
We see that the fixed scheme either consists of two points $c_1,c_2\in C$ with multiplicity  $m=2$,
or a single point $c_1\in C$ of multiplicity $m=4$.
Write $C_i=g^{-1}(c_i)_\red$ for the ensuing reduced fibers. The Canonical Bundle Formula (\cite{Bombieri; Mumford 1977}, Theorem 2)
gives
$ \omega_Y=g^*(\O_{\PP^1}(d-2))\otimes \O_Y(\sum a_iC_i)$ with certain coefficients $0\leq a_i\leq m-1$.
Using $\omega_Y=\O_Y$  we conclude that $\O_Y(\sum a_iC_i)$ is globally generated.
If there is a single multiple fiber, we must have $a_1=0$.
If there are two multiple fibers, the coefficients vanish as well: otherwise  $a_1=a_2=1$ and thus $d=1$,
thus $\O_Y(C_1-C_2)=\O_Y$, contradicting $h^0(\O_Y(C_i))=1$. 
Summing up, in both cases we have $a_i=0$ and $d=2$.
Applying the Canonical Bundle Formula again, we see that the torsion part $\shF\subset R^1g_*(\O_Y)$
has length $h^0(\shF)=2$.

From this information we may compute  $h^1(\O_Y(n))$ for any integer $n$: 
The Projection Formula yields
$R^1g_*(\O_Y(n))=\shF\oplus\O_{\PP^1}(n-2)$.
The Leray--Serre spectral  sequence for the elliptic fibration $g:Y\ra \PP^1$ induces 
an exact sequence
$$
0\lra H^1(\PP^1,  \O_{\PP^1}(n-2))\lra H^1(Y,\O_Y(n))\lra H^0(\PP^1,\shF\oplus\O_{\PP^1}(n-2))\lra 0.
$$
This gives $h^1(\O_Y(2))= 0+3$, whereas $h^1(\O_Y(-2))= 3+2$
Seeking a contradiction, we now suppose that the extension \eqref{bielliptic extension} splits, 
such that  $\Omega^1_Y=\O_Y(2)\oplus\O_Y(-2)$.
This gives $6=h^1(\Omega^1_Y)=h^1(\O_Y(2))+h^1(\O_Y(-2))=8$, contradiction.
\qed
 
\medskip
The extension class for \eqref{bielliptic extension} lies in $\Ext^1(\O_Y(-2),\O_Y(2))=H^1(Y,\O_Y(4))$,
which has dimension $h^1(\O_Y(4))=5$.  It would be interesting to describe this extension class 
explicitly. Note also that the values $h^i(\Theta_Y)$ for $Y=(E\times C)/G$ where
$C$ is elliptic and the $G$-action on it has a fixed point where computed by Partsch
(\cite{Partsch 2013}, Proposition 6.1). It would be interesting to understand the situation in families.

We now   apply our general results on proper group schemes:

\begin{theorem}
\mylabel{bielliptic no lift}
Suppose that $p=2$, and that the bielliptic surface $Y=(E\times C)/G$ is formed
with the rational cuspidal curve $C$ and the group scheme $G=\alpha_2$.
Then $Y$ does not lift to the ring of Witt vectors $W$.
\end{theorem}

\proof
According to \cite{Bombieri; Mumford 1977}, discussion on page 25,
the group scheme $P=\Pic^0_{Y/k}$ has dimension one and embedding dimension two.
The fibration $Y\ra B=E/G$ gives an inclusion $B\subset P$, and the 
$G$-torsor  $X\ra Y$ yields an inclusion $\alpha_2\subset P$.
The presence of multiple fibers shows that $\alpha_2\cap B=0$, and we conclude
that the resulting inclusion $B\times\alpha_2\subset P$ is an equality. In turn, $P/P_\red=\alpha_2$.
Moreover, we have $b_1=2$,  $h^1(\O_Y)=2$ and $h^1(\O_Y)=1$. 
Thus Theorem \ref{picard no lift} applies, and we see that the scheme $Y$ does not lift to the ring $W$.
\qed

%===========================================================
%\appendix
\section{Some homological algebra}
\mylabel{Homological algebra}

In this final section we discuss the relevant homological algebra used throughout the paper.
Our goal is to give a concise description how  
splittings in the derived category, Yoneda extensions and  certain  diagrams   are 
related.  The material should be of independent interest.

Let  $f:M\ra N$ be a homomorphism between two objects $M,N$ in some abelian category $\calA$.
We may regard it as a \emph{two-term complex}, with $f$ as differential. Let $H^0=\Kernel(f)$ and
$H^1=\Cokernel(f)$ be its cohomology, and write  $B=\Image(f)$ for the coboundaries.
Now let $E$ be another object, and $M\stackrel{h}{\ra}E\stackrel{g}{\ra} N$ be some homomorphisms.
We say that $(E,g,h)$ is a \emph{diagram completion} if the diagram
\begin{equation}
\label{cartesian and cocartesian}
\begin{CD}
M	@>\pr>>	B\\
@VhVV		@VViV\\
E	@>>g>	N
\end{CD}
\end{equation}
is \emph{cartesian and cocartesian}. Here $\pr:M\ra B$ and $i:B\ra N$ are the canonical projections and injections, respectively.
The   condition means that  $f=g\circ h$, and that the sequence
$0\ra M\stackrel{(h,\pr)}{\ra} E\oplus B\stackrel{(g,-i)}\ra N\ra 0$
is exact. It follows that
$h:M\ra E$ is a monomorphism, $g:E\ra N$ is an epimorphism, and we have identifications
\begin{equation}
\label{identifications}
\Kernel(g)=\Kernel(\pr)=H^0\quadand \Cokernel(h)=\Cokernel(i)=H^1,
\end{equation}
according to \cite{Kashiwara; Schapira 2006}, Lemma  8.3.11.
The composition of the inclusion $H^0\subset M$ with $h:M\ra E$ yields an inclusion $H^0\subset E$.
In turn, we get a   diagram
\begin{equation}
\label{diagram quasi-isomorphisms}
\begin{CD}
M	@>h>>		E		@<\can<<		H^0\\
@VfVV			@V(g,0)VV				@VV0V\\
N	@>>(\id_N,0)>	N\oplus H^1	@<<(0,\id)<	H^1,
\end{CD}
\end{equation}
and one easily checks that it is commutative. We now regard the vertical maps as two-term complexes,
and the horizontal maps as morphisms between complexes.
Using the identifications \eqref{identifications}, we infer that these are \emph{quasi-isomorphisms}.
We thus may regard \eqref{diagram quasi-isomorphisms} as an \emph{isomorphism}
$$
H^0\oplus H^1[-1]\, \lra \,(M\stackrel{f}{\ra} N)
$$
\emph{in the derived category} $D^b(\calA)$. Note that this constitutes  a  \emph{splitting} of the  complex $M\stackrel{f}{\ra}N$ 
in the sense of Deligne and Illusie (\cite{Deligne; Illusie 1987}, Section 3). 

Recall that by Yoneda's construction \cite{Yoneda 1954}, the groups $\Ext^n(A,B)$ can be defined via equivalence classes
of exact sequences $0\ra C_{n+1}\ra\ldots\ra C_0\ra 0$
with $C_0=A$ and $C_{n+1}=B$.
This   works without the existence of injective
or projective resolutions, and    yields  a $\partial$-functors in $B$.
For details we refer to \cite{Mitchell 1965}, Chapter VII. 
Write $\cl(C_\bullet)\in\Ext^n(A,B)$ for the resulting \emph{Yoneda class}. 
In particular, the horizontal exact sequence in the commutative diagram
\begin{equation}
\label{yoneda ext}
\begin{tikzcd}[row sep = tiny, column sep = small]
		&			& 0\arrow[rd]				&		& 0\\
		&			&					& B\arrow[rd,"i"]\arrow[ru]\\
0\arrow[r]	& H^0\arrow[r]		& M\arrow[rr,"f"']\arrow[ru,"\pr"]	&		& N\arrow[r]	& H^1\arrow[r]	& 0\\
\end{tikzcd}
\end{equation}
yields a Yoneda class, which we denote by $\cl(f)\in \Ext^2(H^1,H^0)$.
It coincides with the \emph{Yoneda product} $\cl(M)\ast\cl(N)$ of the   
\emph{extension classes} for the two short exact sequences with kinks.

\begin{lemma}
\mylabel{completion yoneda class}
The homomorphism $f:M\ra N$ admits a diagram completion $(E,g,h)$ if and only if the Yoneda class
$\cl(f)\in\Ext^2(H^1,H^0)$ vanishes.
\end{lemma}

\proof
In somewhat different formulation, this already appears in \cite{Barakat; Bremer 2008}, Theorem 5.1. Let me give an independent argument.
The short exact sequence to the right in \eqref{yoneda ext} yields
an extension class $\cl(N)\in\Ext^1(H^1,B)$, whereas the short exact sequence to the left gives
a long exact sequence
\begin{equation}
\label{ext sequence}
\Ext^1(H^1,M)\lra \Ext^1(H^1,B)\stackrel{\partial}{\lra} \Ext^2(H^1,H^0).
\end{equation}
By definition of this sequence (\cite{Mitchell 1965}, Chapter VII, Section 5), the image of the extension class $\cl(N)$ under the connecting map 
is   $\cl(f)\in\Ext^2(H^1,H^0)$.

Suppose $\cl(f)=0$. Then the extension
$0\ra B\ra N\ra H^1\ra 0$ arises from an extension $0\ra M\stackrel{h}{\ra} E\ra H^1\ra 0$,
and this means that there is a cocartesian diagram 
\begin{equation}
\label{yoneda and completions}
\begin{CD}
M	@>h>>	E\\
@V\pr VV		@VVgV\\
B	@>>i>	N.
\end{CD}
\end{equation}
Using that $h:M\ra E$ is a monomorphism, together with \cite{Kashiwara; Schapira 2006}, Lemma 8.3.11
we infer that the above cocartesian diagram is also cartesian. In turn, $(E,g,h)$ is a diagram completion.

Conversely, suppose there is a diagram completion $(E,g,h)$, giving a cartesian and cocartesian diagram 
\eqref{yoneda and completions}.
Now recall  that  $h$ is a monomorphism and   $\Cokernel(h) = \Cokernel(i)=H^1$.
This means that the extension class $\cl(N)$ lies in the image of the   map on the left in \eqref{ext sequence},
and thus   $\cl(f)\in\Ext^2(H^1,H^0)$ vanishes.
\qed

The diagram completions form a category $\operatorname{Cp}(M\stackrel{f}{\ra}N)$.
In this category, the morphisms $ (E,g,h)\ra (E',g',h')$ are those  homomorphisms $\varphi:E\ra E'$ making the diagram
\begin{equation}
\label{homomorphism diagram}
\begin{tikzcd}[row sep=tiny ]
		&							& E\oplus B\arrow[rd,"{(g,-i)}"]\ar[dd,"\varphi\oplus\id_B"]\\		
0\arrow[r]	& M\arrow[ru,"{(h,\pr)}"]\ar[rd, "{(h',\pr')}"']	& 			& N\ar[r]		& 0\\
		&							& E'\oplus B\arrow[ru,"{(g',-i')}"']	
\end{tikzcd}
\end{equation}
commutative. In other words, $\varphi\oplus\id_B$ is a morphism of extensions.
The latter is an isomorphism, by the Five Lemma, hence the same holds for the direct summand $\varphi$,
so the category $\Cp(M\stackrel{f}{\ra}N)$ is a \emph{groupoid}.

Each homomorphism $\xi:H^1\ra H^0$ yields an endomorphism
\begin{equation}
\label{epi mono}
E\stackrel{g}{\lra} N \stackrel{\can}{\lra} H^1\stackrel{\xi}{\lra} H^0\stackrel{\can}{\lra} M\stackrel{h}{\lra} E,
\end{equation}
which we denote by $\xi_E:E\ra E$. Consider the   endomorphism $\xi_E\circ\xi_E$. By definition, it
is a composition of the form $\ldots\ra H^0 \stackrel{\can}{\ra} M \stackrel{h}{\ra} E \stackrel{g}{\ra} N\ra\ldots$.
This vanishes, because the composition of $f=g\circ h$ with the inclusion
of $H^0=\Kernel(f)\subset M$ is zero.  
It follows that $\id_E+\xi_E$ is an automorphism of $E$, with inverse $\id_E-\xi_E$.
From the definition and \eqref{homomorphism diagram} one easily infers that it is actually an automorphism of $(E,g,h)$.
In turn, we obtain an   homomorphism of groups 
\begin{equation}
\label{automorphism groups}
\Hom_{\calA}(H^1,H^0)\lra \Aut(E,g,h)\subset \Aut_\calA(E),\quad \xi\longmapsto \id_E+\xi_E
\end{equation}
This map is injective, because in \eqref{epi mono}, the two arrows to the left are epimorphisms,
whereas the two arrows to the right are monomorphisms.

\begin{proposition}
\mylabel{automorphisms equal}
The inclusion $\Hom_{\calA}(H^1,H^0)\subset \Aut(E,g,h)$ is an equality   provided that the canonical inclusion
$H^0\subset M$ admits a retraction.
\end{proposition}

\proof 
Choose a retraction $r:M\ra H^0$. Using the functoriality of our maps, we may assume
$M=H^0\oplus B$. The composition $h\circ r:B\ra E$, together with the 
universal property of cocartesian diagrams shows that the surjection $E\ra N$ admits a splitting. Thus we may assume
$E=H^0\oplus N$, where the morphism $h:M\ra E$ is given by the matrix $(\begin{smallmatrix} \id&0\\0&\pr\end{smallmatrix})$,
and $g:E\ra N$ is given by $(0,\id)$. Using \eqref{homomorphism diagram}, one sees that each automorphism of 
$(E,g,h)$ is of the form $\id_E+\xi_E$ for some homomorphism $\xi:H^1\ra H^0$.
\qed

\medskip
For the dual situation, we make the following observation:

\begin{proposition}
\mylabel{objects isomorphic}
The category $\Cp(M\stackrel{f}{\ra}N)$ has precisely one isomorphism class
provided that the canonical surjection $N\ra H^1$ admits a section.
\end{proposition}

\proof
Fix a section $s:H^1\ra N$, and write  $N=B\oplus H^1$. Set $E_0=M\oplus H^1$.
Let $g_0:M\ra E_0$ be the canonical inclusion and $h_0:E=M\oplus H^1\ra B\oplus H^1=N$ by the matrix $(\begin{smallmatrix}\pr&0\\0&\id\end{smallmatrix})$.
One easily checks that $(E_0,g_0,h_0)$ is a diagram completion.
Let $(E,g,h)$ be another diagram completion. Composing $g:E\ra N$ with the retraction $N\ra B$ and using
the universal property of cartesian squares, we get $E\simeq M\oplus H^1$, and infer that $(E,g,h)$ is isomorphic to $(E_0,g_0,h_0)$.
\qed

\medskip
Now we bring in topology.
Suppose that $\calC$ is a ringed site. For the sake of exposition, we assume that there is a final object $X\in\calC$, and
write the structure sheaf as $\O_X$. We regard
the objects as ``open sets'', and write them as $U\ra X$.
From now on we assume  that our abelian category is $\calA=(\O_X\text{-Mod})$, such that our 
$f:M\ra N$ is a homomorphism of $\O_X$-modules. Note also that there are enough injective objects. 
\emph{Furthermore assume that $H^1=\Cokernel(f)$ 
is locally free of finite rank, and that $\pr:M\ra B$ locally admits sections,} as in \cite{Deligne; Illusie 1987}, Section 3.2.

The first condition ensures  that  the contravariant functor 
$U\mapsto \Hom_{\O_U}(H^1|U,F|U)$ satisfies the sheaf axiom, where  $F$ is an abelian sheaf, 
and $U\ra X$ runs over the objects of $\calC$. We denote the resulting sheaf $\uHom_{\O_X}(H^1,F)$.
Note also that we  have an identification
$\Ext^n(H^1,F)= H^n(X,\uHom_{\O_X}(H^1,F))$, $ n\geq 0$
of universal $\partial$-functors in $F\in\calA$.

Let $\shCp(M\stackrel{f}{\ra}N)$ be the \emph{category fibered in groupoids} over $\calC$, whose 
objects over $U\ra X$ are    the diagram completions $(E,g,h)$ 
for the restrictions $M|U\stackrel{f|U}{\ra} N|U$.
Morphisms $(E,g,h)\ra (E',g',h') $   over a given $U\ra U'$ are isomorphisms 
$(E,g,h)\ra (E'|U,g'|U,h'|U) $. One easily checks that this category is fibered and satisfies the stack axioms.
Roughly speaking, this means that all Hom presheaves are sheaves, and that all descend data are effective.
See \cite{Olsson 2016}, Chapter 2 and 3 for the relevant definitions.

For each object $(E,g,h)$ over $U\ra X$, we obtain from \eqref{automorphism groups} a homomorphism of group-valued sheaves
$$
\Psi_{E,g,h}:\uHom_{\O_X}(H^1,H^0)|U\lra \Aut_{(E,g,h)/U},\quad \xi\mapsto \id_E+\xi_E
$$
We observe:

\begin{proposition}
\mylabel{gerbe of completions}
The above  are isomorphisms, and the stack $\shCp(M\stackrel{f}{\ra}N)$ is a gerbe banded by  the abelian sheaf 
$\uHom_{\O_X}(H^1,H^0)$.
\end{proposition}
 
\proof
We have to check that all $\Psi_{E,g,h}$ are      isomorphisms, and that all objects in $\shCp(M\stackrel{f}{\ra}N)$
are locally isomorphic. Both are local problems, and by our overall assumptions it suffices to treat the case that
$N\ra H^1$ and $M\ra B$ admit splittings. The assertion on the homomorphisms and the objects follow from 
Proposition \ref{automorphisms equal} and \ref{objects isomorphic}, respectively.
\qed

\medskip
Recall that Deligne and Illusie (\cite{Deligne; Illusie 1987}, Section 3) defined the \emph{gerbe of splittings}, 
which we denote by $\shSc(M\stackrel{f}{\ra}N)$.
The objects over $U\ra X$ are the splittings $s$ for the canonical projection $N|U\ra H^1|U$,
and the morphisms $s\ra s'$  between two splittings are defined as the homomorphisms $\xi:H^1|U\ra M|U$
with $s'= s + \pr\circ \xi$. Via the tautological map
$$
\Phi_{E,g,h}:\uHom_{\O_X}(H^1,H^0)|U\lra \Aut_{s/U},\quad \xi\longmapsto \xi
$$
this also becomes a gerbe banded by $\uHom_{\O_X}(H^1,H^0)$.
In turn, we have two gerbes banded by the same coefficient sheaf, giving two cohomology classes. Our main result here is:

\begin{theorem}
\mylabel{same gerbe}
The gerbe $\shCp(M\stackrel{f}{\ra}N)$ of diagram completions and the gerbe 
$\shSc(M\stackrel{f}{\ra}N)$ of splittings have the same class in 
the cohomology group
$$
H^2(X,\uHom_{\O_X}(H^1,H^0))=\Ext^2(H^1,H^0).
$$
Moreover, either  of them admits a global object if and only if the Yoneda class $\cl(f)\in \Ext^2(H^1,H^0)$
of the exact sequence
$0\ra H^0\ra M\stackrel{f}{\ra} N\ra H^1\ra 0$   vanishes.
\end{theorem}

\proof
First we construct a functor from the latter category to the former.
Suppose we have a global splitting $s:H^1\ra N$, and write $N=B\oplus s(H^1)$, where $s(H^1)=\Image(s)$. Set $E=M\oplus s(H^1)$.
Let the homomorphisms $g,h$ be defined by the   diagram
\begin{equation}
\label{functor}
\begin{CD}
M		@>\pr>>		B\\
@V(\id,0)VV			@VV(\id,0)V\\
M\oplus s(H^1)	@>>(\pr,\id)>	B\oplus s(H^1).
\end{CD}
\end{equation}
Clearly, this constitutes a diagram completion.
The same reasoning applies locally over   $U\ra X$. According to \cite{Giraud 1971}, Chapter IV, Corollary 2.2.7,
the functor $s\mapsto (E,g,h)$ is an equivalence of categories.
In turn, if one of them admits a global object, so does the other.
By Proposition \ref{completion yoneda class}, the category of diagram completions contains a global object
if and only if the Yoneda class vanishes.

It remains to check that the gerbe classes coincide, and do not differ by a sign, say.
For this we have to check that they are banded by
the coefficient sheaf $\uHom_{\O_X}(H^1,H^0)$ in the same way.
Our   construction \eqref{functor} is functorial in $s$. In particular, each homomorphism $\xi:H^1\ra H^0$, viewed as an automorphism
of $s$, yields the automorphism of $(E,g,h)$ with $E=M\oplus s(H^1)$ as 
above given by the matrix $(\begin{smallmatrix}\id&\xi\\0&\id\end{smallmatrix})=\id_E+\xi_E$.
In light of \eqref{automorphism groups}, the  actions of the abelian sheaf $\uHom_{\O_X}(H^1,H^0)$ via $\Psi$ and $\Phi$
on the objects $s$ and $(E,g,h)$ coincide.  
\qed

\medskip
Let us close with the following remark: Since the quasi-isomorphisms in the homotopy category of cochain complexes
admit a calculus of left and right fractions in the sense of Gabriel and Zisman \cite{Gabriel; Zisman 1967},
any   isomorphism  in the derived category represented by quasi-isomorphisms as in the diagram \eqref{diagram quasi-isomorphisms}
may also be represented by 
quasi-isomorphisms $H^0\oplus H^1[-1]\leftarrow C^\bullet\ra (M\stackrel{f}{\ra}N)$, with arrows pointing in 
reverse directions. This dichotomy seems to lie at the heart of the matter for the preceding results.

%===========================================================

\end{document}